\documentclass[12pt, a4paper, reqno, twoside, makeidx]{amsart}
\usepackage[margin=2.8cm, marginpar=2cm]{geometry}
\usepackage{graphicx}
\usepackage{color}
\usepackage{amsmath, amsthm, amssymb, amsfonts, amscd}
\usepackage[english]{babel}
\usepackage[all]{xy}
\usepackage{url}
\usepackage{hyperref}
\usepackage{mathtools}

\theoremstyle{remark}
\newtheorem{rmk}{Remark}
\theoremstyle{plain}
\newtheorem{thm}{Theorem}

\newtheorem{defi}[thm]{Definition}

\theoremstyle{definition}
\newtheorem{exa}[thm]{Example}
\newcommand{\N}{\mathbb{N}}

\newcommand{\R}{\mathbb{R}}

\newcommand{\HH}{\mathrm{H}}

\author{Daniele Celoria}
\author{Barbara I. Mahler}

\title{A statistical approach to knot confinement via persistent homology}
\begin{document}
\maketitle

\begin{abstract}
In this paper we study how randomly generated knots occupy a volume of space using topological methods. To this end, we consider the evolution of the first homology of an immersed metric neighbourhood of a knot's embedding for growing radii. Specifically, we extract features from the persistent homology of the Vietoris-Rips complexes built from point clouds associated to knots. Statistical analysis of our data shows the existence of increasing correlations between geometric quantities associated to the embedding and persistent homology based features, as a function of the knots' lengths. We further study the variation of these correlations for different knot types. Finally, this framework also allows us to define a simple notion of deviation from ideal configurations of knots. 
\end{abstract}

\section{Introduction}\label{sec:intro}

In this paper we consider the following question: does the topology of a random knot, \emph{i.e.}~its knot type, influence how it ``occupies'' space? More specifically, do more complicated knots tend to be more compact or loose with respect to their simpler counterparts of the same length? Similar ideas have previously been considered in \emph{e.g.}~\cite{polyconfinement,polyconfinement2,polyconfinement3,spectrumconfinement}, often with the main goal of understanding the mechanism of DNA-- or, more generally, polymer packing in a confined volume. Prominent instances of such studies are concerned with tight DNA packing in viral capsids~\cite{phages,capsids,capsids2,capsids3} and DNA packing in cells (see \emph{e.g.}~\cite{packing, packing2} and the review~\cite{biopolymer}). We refer to the comprehensive survey~\cite{micheletti} for further references and several related notions.

Our main tool to address these matters is persistent homology (PH) \cite{ghrist2008barcodes, carlsson2009}; this is a relatively new technique in topological data analysis, commonly used to detect insightful topological and geometric features of point clouds. Roughly speaking, PH associates a filtered simplicial complex to a point cloud. The filtered homology groups of these chain complexes often capture subtle properties of the point cloud~\cite{edelsbrunner2010computational}.

We start by generating random piecewise linear knots with prescribed length and\slash or topology. We then create a point cloud for each knot $K$ by linearly interpolating between the endpoints of the PL curve, and we compute the PH of its Vietoris-Rips filtration in dimension $1$. 
Intuitively, we use PH to examine the changes in topology occurring in a metric neighbourhood of a random knot, when the radius of this neighbourhood varies in $\R_{\ge 0}$.  

We extract several features from the obtained PH and quantify the variation in the correlation between these features and either the volume ``occupied'' by the knot or the average crossing number (ACN) for increasing lengths and for all knot types with up to six crossings.  The most prominent feature we extract from PH is $\mathcal{I}(K)$, the integral of the Betti curve of a persistence diagram, which might be of independent interest.

Note that unlike previous approaches, such as those outlined in~\cite[Sec.~8]{micheletti}, we do not prescribe the geometry of the confined volume, but rather work backwards, by first generating the knots and only afterwards analysing their relationship with the minimal volumes that contains them. More precisely, we consider different kinds of measures for the space ``filled in'' by a knot: the volume of the circumscribing sphere, \emph{i.e.}~the volume of the smallest sphere that encloses the knot, the volume of the convex hull determined by the knot, and the radius of gyration. This last quantity is often used as a meaningful and computationally convenient measure of compactness of proteins and polymers (see \emph{e.g.}~\cite{gyrationcompactness,gyrationdiffusivity,gyrationselfavoiding,stamech, micheletti}). We remark that, in the case of proteins, compactness is defined as the ratio of the accessible surface area of a protein to the surface area of the sphere of the same volume. We instead consider a more intuitive and geometric notion of compactness. \newline

Most prominently, we show the existence of an inverse correlation between the integral $\mathcal{I}(K)$ and the various notions of volume occupied by $K$ mentioned above, as well as a direct correlation between $\mathcal{I}(K)$ and the average crossing number. The magnitude of these correlations increases for increasing knot lengths. Furthermore, we find that these correlations appear to differentiate between different knot types. Peculiarly, the intensity of the computed correlations does not appear to be directly related to classical measures of a knot's complexity, such as the minimal crossing number.

To better appreciate these results, especially the resulting subdivision into knot types, we also compute the average Betti curve for each knot type considered, as well as their integrals. We observe an almost perfectly linear relation between the average integrals of the Betti curves and the knot lengths; the same holds for the average maxima of the Betti curves, which are in turn correlated to the average number of shallow angles in the embeddings. 
Indeed, these relations show a clear divide among the considered knot types.\newline

We then turn to the related concept of \emph{ideal knot} embeddings. These are special embeddings of knots, whose study was pioneered by Stasiak~\cite{stasiak1998ideal}. Their geometry is particularly simple, in that they minimise the length of a rope (having unit diameter) that is needed to tie a specific knot. We show how the PH framework developed here can be used to define a simple numerical measure of how ``far away'' a given knot is from an ideal embedding.\newline

We make the concepts mentioned in the introduction rigorous in Section~\ref{sec:basicknot}, and give a basic overview of the techniques we use in Section~\ref{sec:PH}. We then detail how we generated our data in Section~\ref{sec:data}, and present the results in Section~\ref{sec:results}. Finally, in Section~\ref{sec:distanceideal} we define the aforementioned deviation from ideal knot embeddings.

\subsection{Acknowledgements:} The authors wish to thank Oxford's Mathematical Institute for  providing access to computational resources, and Agnese Barbensi for interesting conversations. DC was supported by the  European Research Council (ERC) under the European Union’s Horizon 2020 research and innovation programme (grant agreement No 674978).

\section{Knot theory}\label{sec:basicknot}

We call a \emph{knot $K$} the image of a smooth embedding of $S^1$ in $S^3$, and reserve the notation $\mathcal{K}$ to denote  knot types. We refer to~\cite{rolfsen2003knots} for basic definitions in knot theory (see also \cite[Ch.~8]{ohara}). In what follows, we relax the smoothness condition to allow the approximation of a smooth embedding by equilateral polygons. These curves will be referred to as piecewise linear (PL) knot embeddings.

The \emph{length} $\ell(K)$ of a knot $K$ is the usual Euclidean length of $K$. When considering the length in the PL case, we always assume that all segments that compose a given polygonal knot are of equal length; we take each segment of unit length, so that $\ell(K)$ coincides with the number of edges used.

\begin{figure}[ht]
\centering
\includegraphics[width=6.5cm]{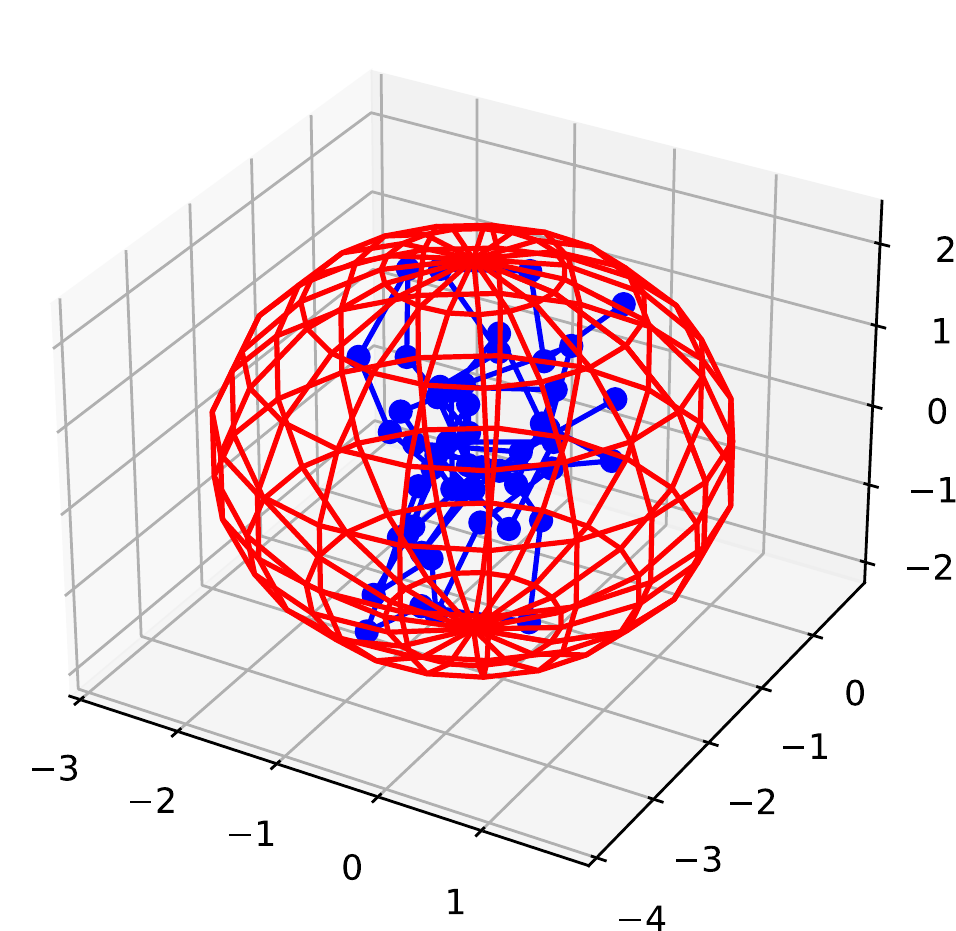}
\includegraphics[width=6.5cm]{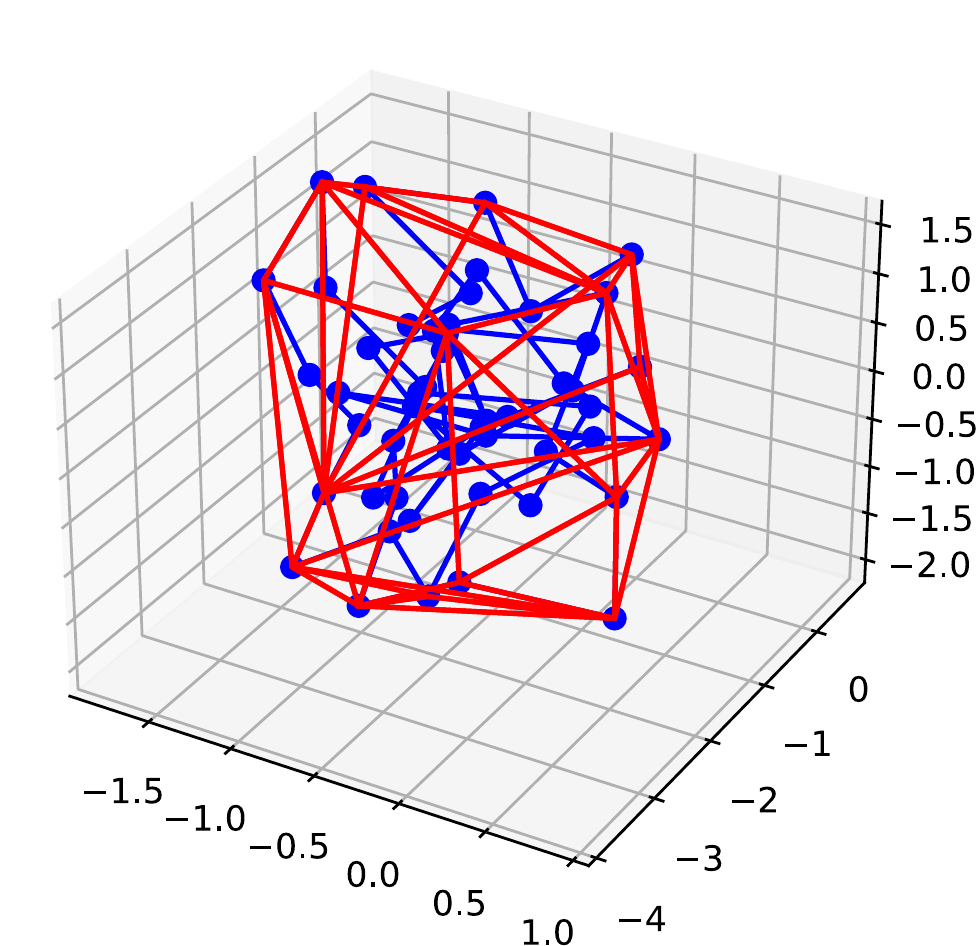}
\caption{In red, a circumscribing sphere and convex hull for the blue PL trefoil knot with $50$ edges. In this example, the sphere has volume $\sim 50$, while the convex hull's volume is $\sim 14$.}
\label{fig:volumes}
\end{figure}

We are interested in investigating how efficiently a given embedding can occupy a volume; we therefore consider different kinds of measures of compactness for a knot embedding. In increasing order of accuracy, for each $K$ we compute the volume of the minimal sphere and the convex hull surrounding $K$ (see Figure~\ref{fig:volumes}). We also take into consideration the radius of gyration $Rg$ of a PL-embedded knot.

We denote by $\nu_t (K)$ the metric neighbourhood of $K$ of radius $t >0$. We obtain this by considering the union of the radius $t$  discs contained in the affine planes $P_p + p$ that are centred at the points $p$ in the image of the embedding, and that are orthogonal to the embedding. 
We crucially point out that in what follows we will not necessarily only consider regular (\emph{i.e.}~non  self-intersecting) neighbourhoods (see Figure~\ref{fig:8unknot}). 

The \emph{injectivity radius} $IR(K)$ of $K$ (see \cite{litherland1999thickness}) is the supremum among all radii $t \ge 0$ for which the tubular neighbourhood $\nu_t(K)$ is regular. In other words, $IR(K)$ is the smallest value of the neighbourhood's radius such that $\nu_t(K)$ comes into contact with itself. Intuitively, if we regard the given knot $K$ as being made of rope, $IR(K)$ represents the maximal radius that the rope can have while being knotted in the ``shape'' of $K$.

The \emph{length over diameter ratio} of $K$ is the quotient  $$L\slash D (K) = \frac{\ell(K)}{2IR(K)}$$ between the length of $K$ and  twice the injectivity radius. The ratio $L\slash D$ was  introduced and has been extensively studied by Stasiak~\cite{stasiak1998ideal} (the inverse of this quantity is also known as the \emph{thickness} of $K$~\cite{litherland1999thickness}), due to its relation with the notion --- also introduced by Stasiak --- of ideal knot.

An \emph{ideal knot} (see the monograph \cite{stasiak1998ideal}) is an embedded knot that minimises the $L\slash D$ ratio within its knot type. Note that a priori there might be more than one ideal representative of a given knot type.  One interpretation of $L\slash D$ for an ideal knot $K$ is as the minimal length needed to tie a knot of the knot type of $K$ with a rope that has diameter $1$.

We will also use other quantities that can be associated with PL knot embeddings: a discrete analogue of torsion and curvature, and the \emph{average crossing number}~(ACN). This latter quantity is defined as the integral of the function $S^2 \longrightarrow \N$ associating to a point $p$ on the unit sphere the crossing number of the diagram obtained by projecting $K$ onto the plane tangent to the sphere at $p$ (strictly speaking, we should also renormalise by dividing by $4\pi$ and restrict to projections with at most double points as singularities). \newline

We have just seen that the notion of injectivity radius is crucial for the definition of several knot properties. One of the key technical aspects of this paper is the use of a straightforward generalisation of $IR$, and that PH can be effectively used to compute it. Given a knot $K$, consider the neighbourhood $\nu_{t}(K)$ for $ t \in [0 , \infty[$. For small enough $t$, the homology of the embedded neighbourhood is of rank $1$ in dimension $1$. The topology changes as soon as we get to $t = IR(K)$, where (generically) the rank of the homology of $\nu_{t}(K)$ increases by one. Similarly, for increasing $t$, we can keep track of all the times $t$ where the topology of $\nu_{t}(K)$ changes. Note that for values of $t$ greater than $RS(K)$, the radius of the circumscribing sphere,  $\mathrm{H}_1(\nu_t(K))$ vanishes, and the only nontrivial homology is of rank $1$ in degree $0$.

We can now introduce the Betti curve for the first homology, which is one of the main objects we will consider in what follows. 
\begin{defi}\label{def:betticurve}
Call the (first) \emph{Betti curve} of a knot $K$ the integer valued function $$\beta_1(K) : \R_{\ge 0} \longrightarrow \N,$$  
defined as $ t \mapsto \mathrm{rk}(\HH_1(\nu_t(K))).$
\end{defi}
It follows from the previous discussion that $\beta_1(K)$ is $1$ for small values of $t$ and becomes definitely $0$ for $t\gg0$.

An interesting property of Betti curves is that (just as with persistent landscapes~\cite{landscapes}) we can add them and take averages; we will take advantage of this fact in Section~\ref{sec:results}.

Note that if $K$ presents some small-scale configuration (such as the smaller twirl in Figure~\ref{fig:8unknot}), then, after increasing the radius more than a certain threshold, its contribution to $\mathrm{rk}(\HH_1(\nu_t(K)))$ vanishes  (after the neighbourhood engulfs the small-scale configuration). 

We will argue shortly that PH can be used to closely approximate $\beta_1(K)$. In fact, by considering distributions of points closely approximating the embedding $K$ and increasing in density, we get increasingly  better approximations of $\beta_1$.

\begin{exa}
Consider the planar standard embedding $\bigcirc$ of $S^1$ in $\R^3$, with radius $R$. Then $\beta_1(\bigcirc) \equiv 1$ on $[0,R[$, and $0$ on $[R, \infty[$. For the rather simple unknotted embedding in Figure~\ref{fig:8unknot}, the Betti curve resembles that of Figure~\ref{fig:betti_of_infty}.
\begin{figure}[ht]
\centering
\includegraphics[width=3cm]{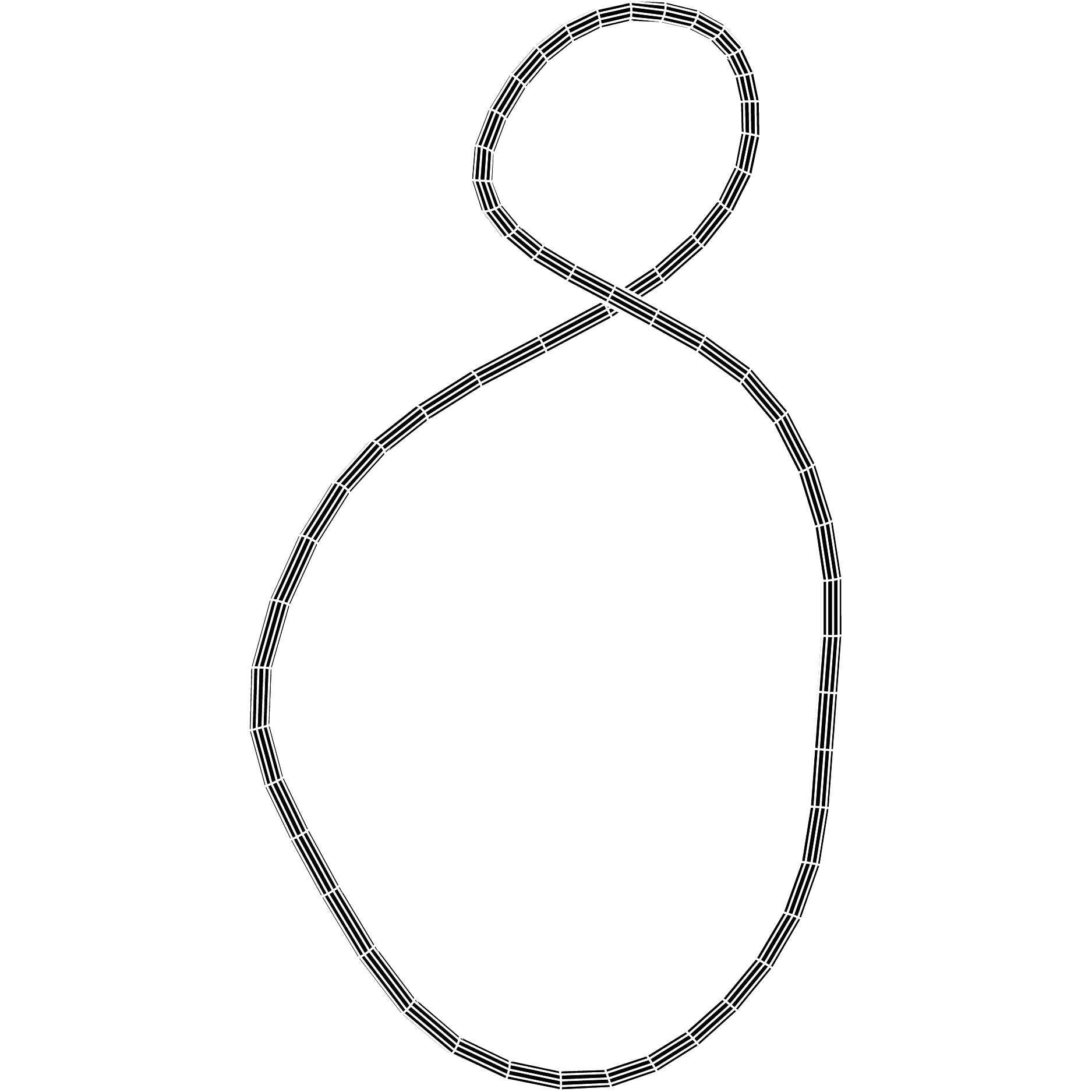}
\includegraphics[width=3cm]{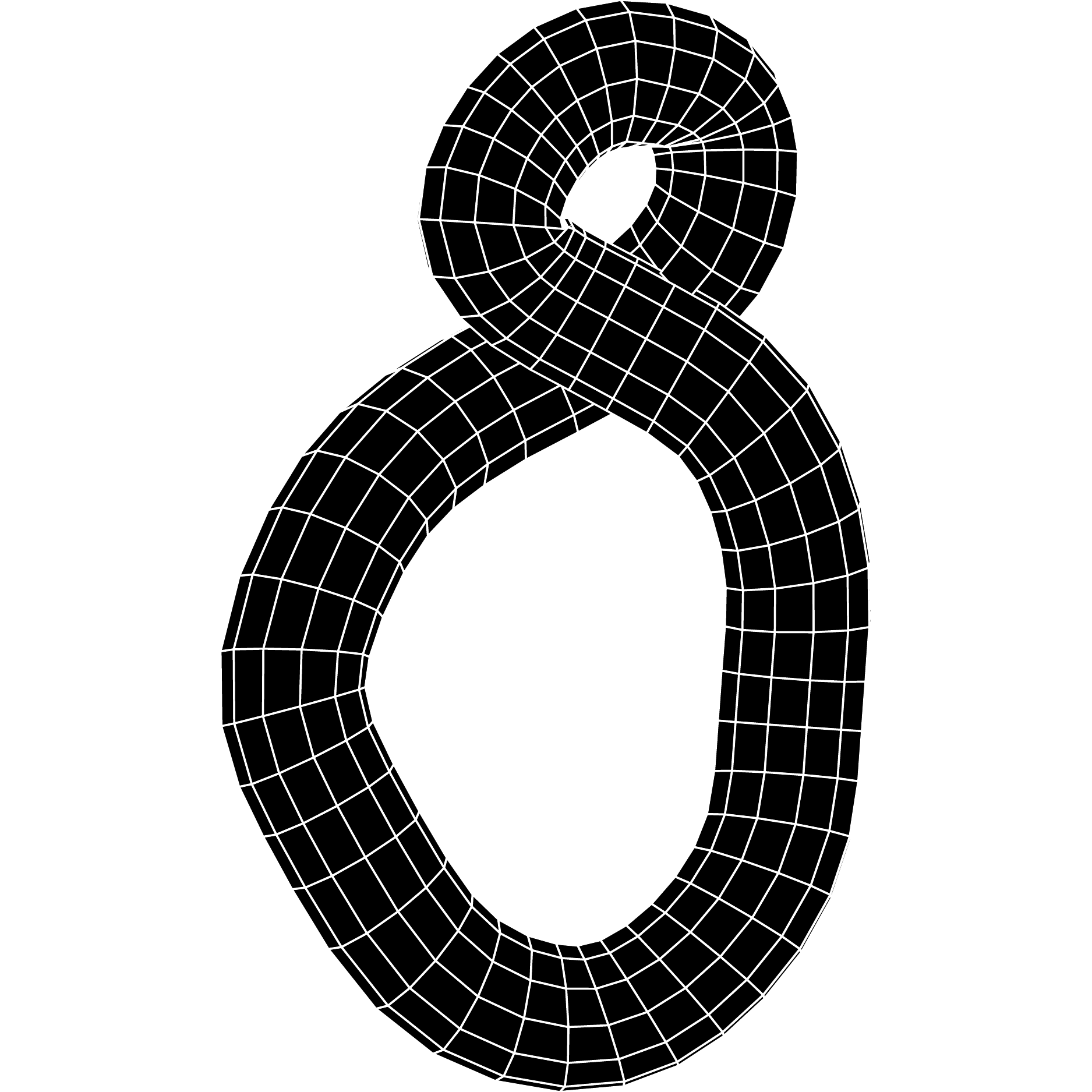}
\includegraphics[width=3cm]{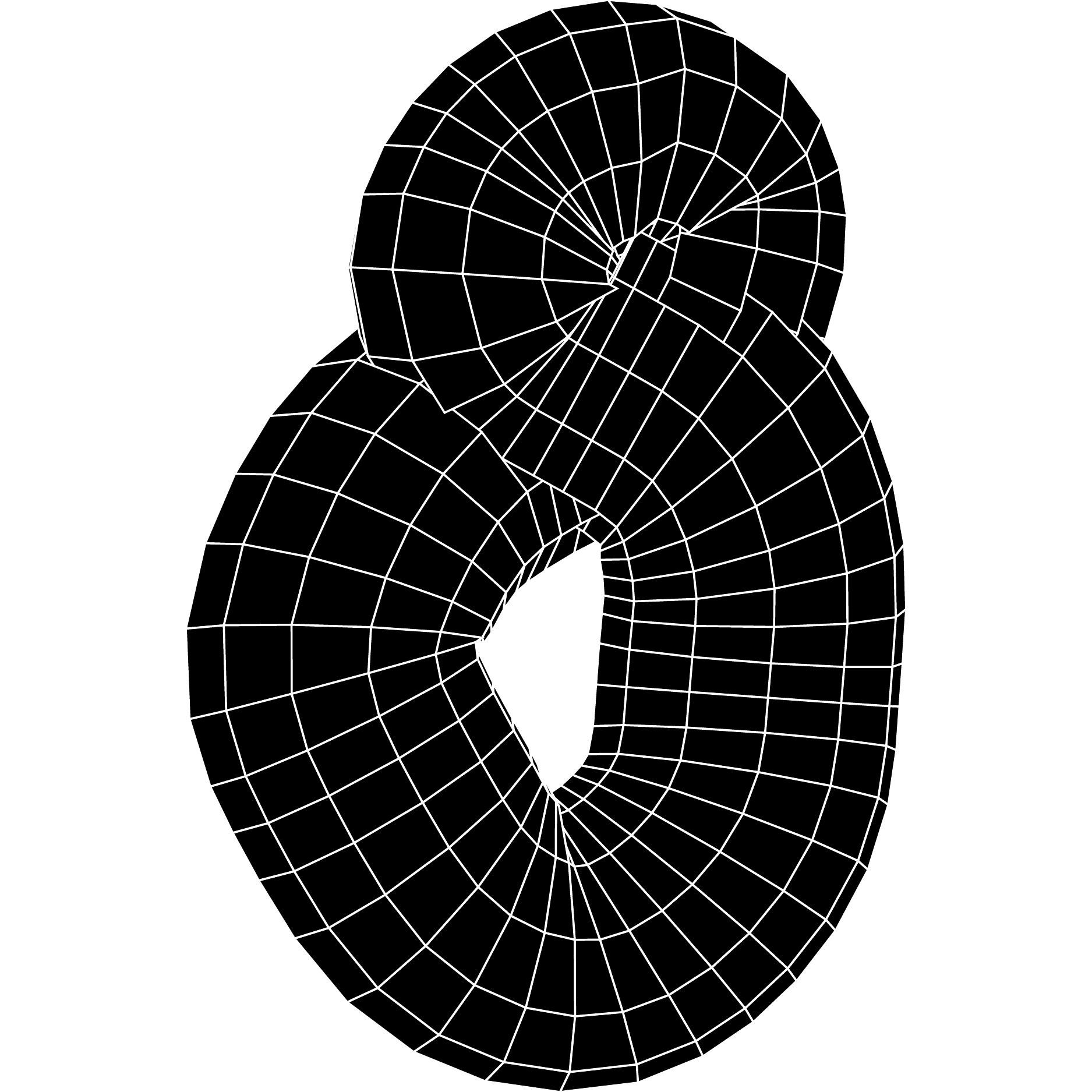}
\includegraphics[width=3cm]{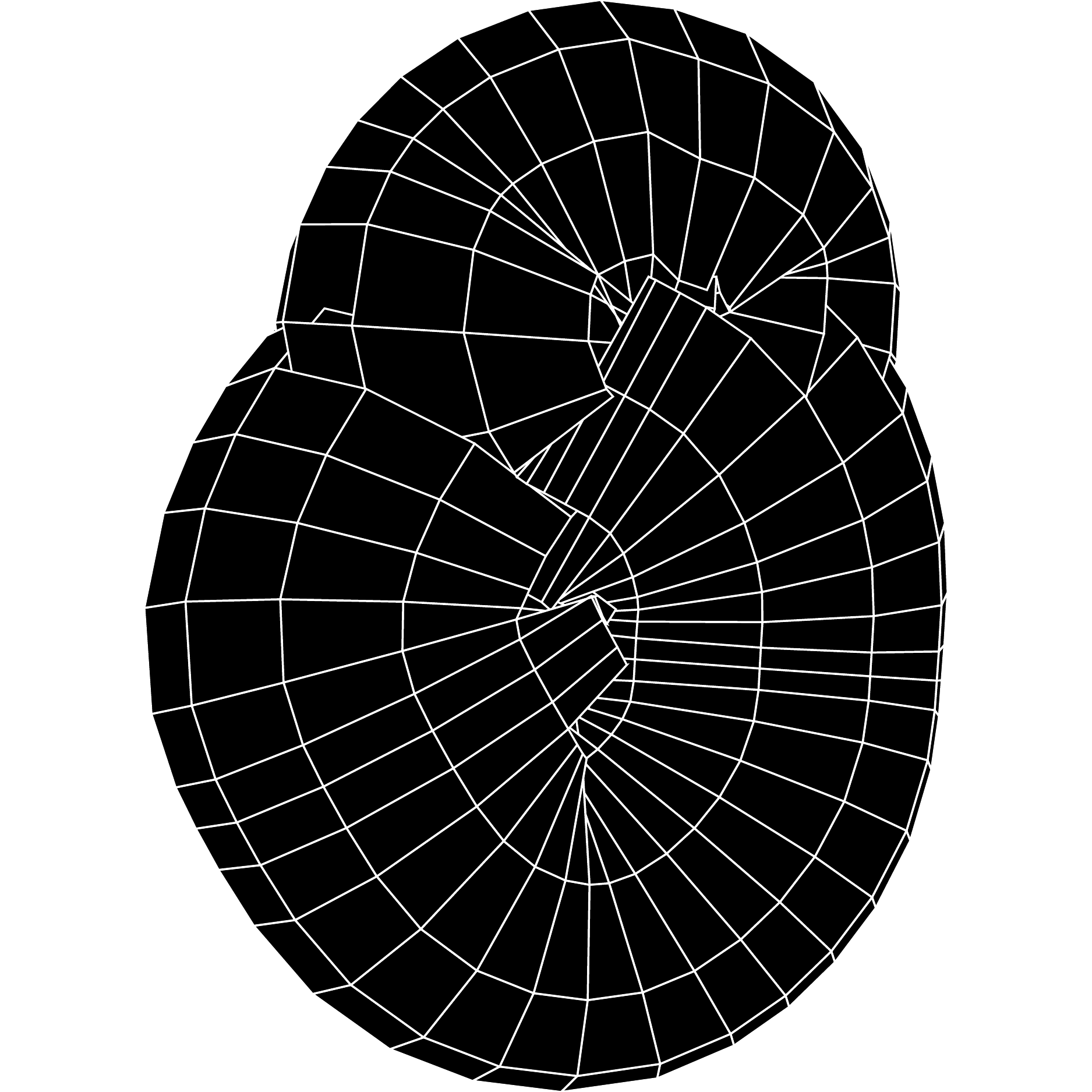}
\caption{The evolution of the neighbourhood $\nu_t(K)$ for an eight-shaped unknot. The last stage on the right is homeomorphic to a $3$-ball. The image was generated using KnotPlot~\cite{knotplot}.}
\label{fig:8unknot}
\end{figure}
\begin{figure}[ht]
\centering
\includegraphics{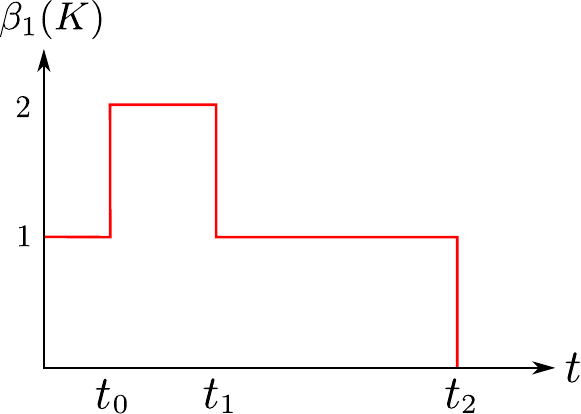}
\caption{A schematic representation of the function $\beta_1(K)$ for the unknotted embedding from Figure~\ref{fig:8unknot}. The values $t_i$ mark the times of $t$ where the topology of $\nu_t(K)$ changes.}
\label{fig:betti_of_infty}
\end{figure}

\begin{figure}[ht]
\centering
\includegraphics[width=3.5cm]{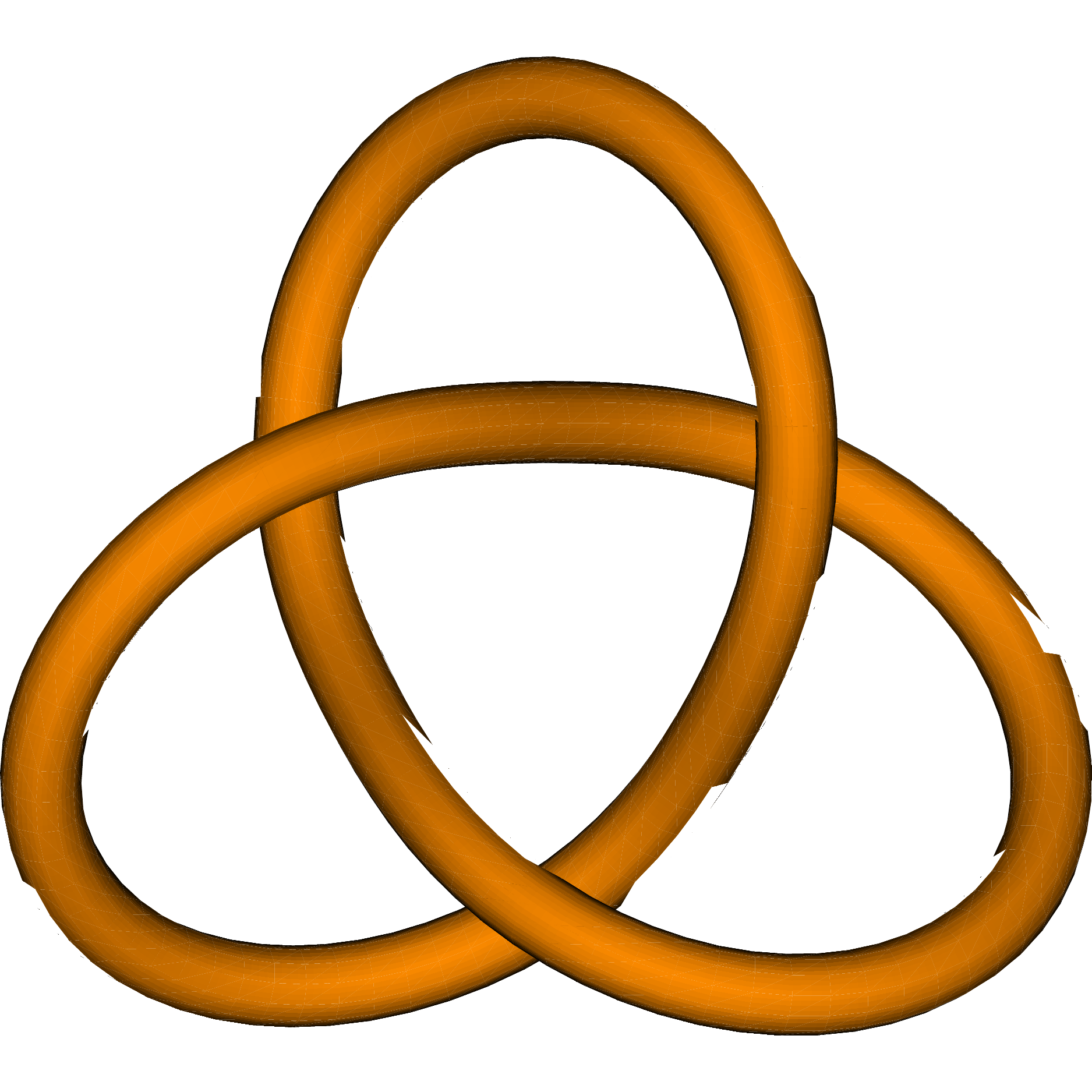}
\includegraphics[width=3.5cm]{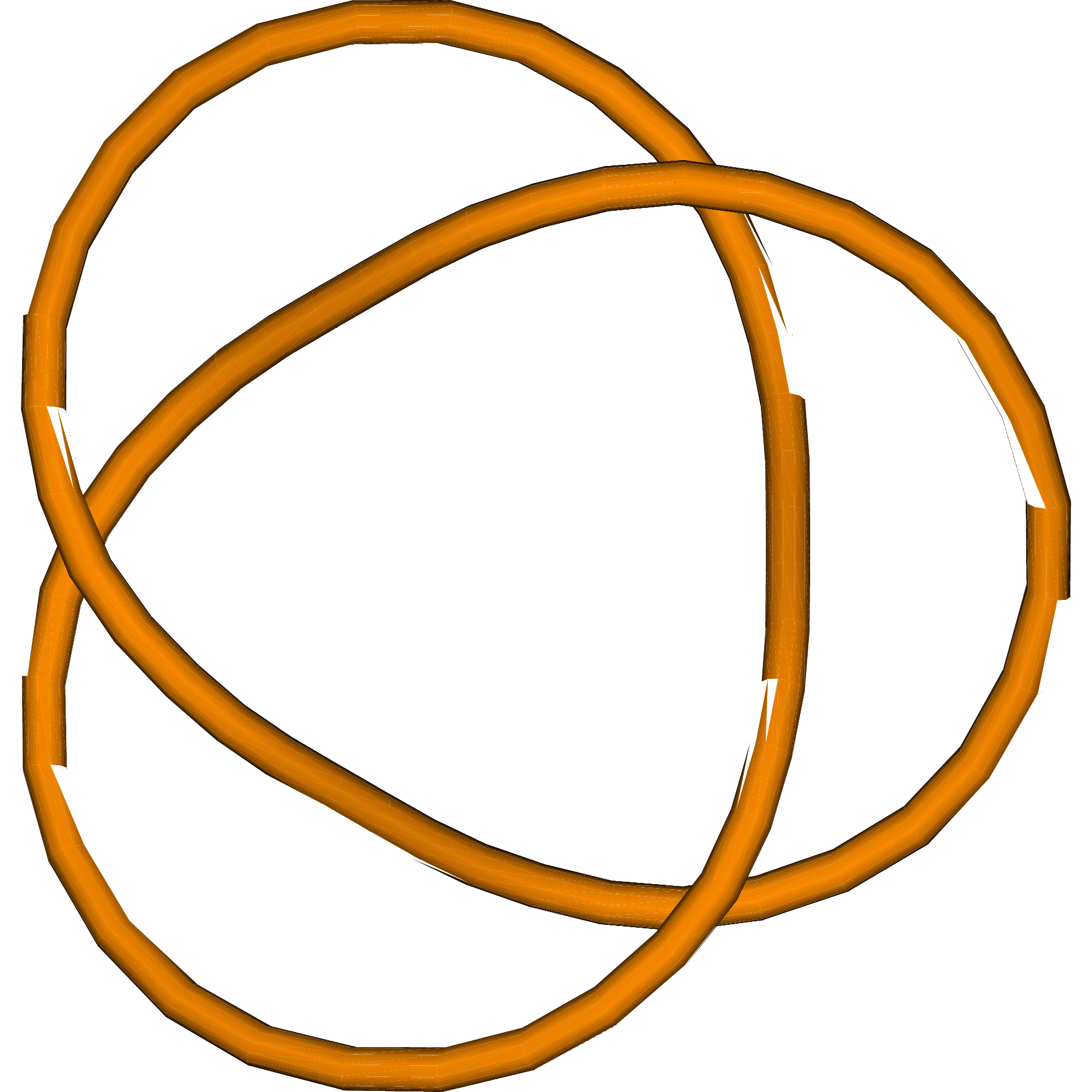}
\includegraphics[width=3.5cm]{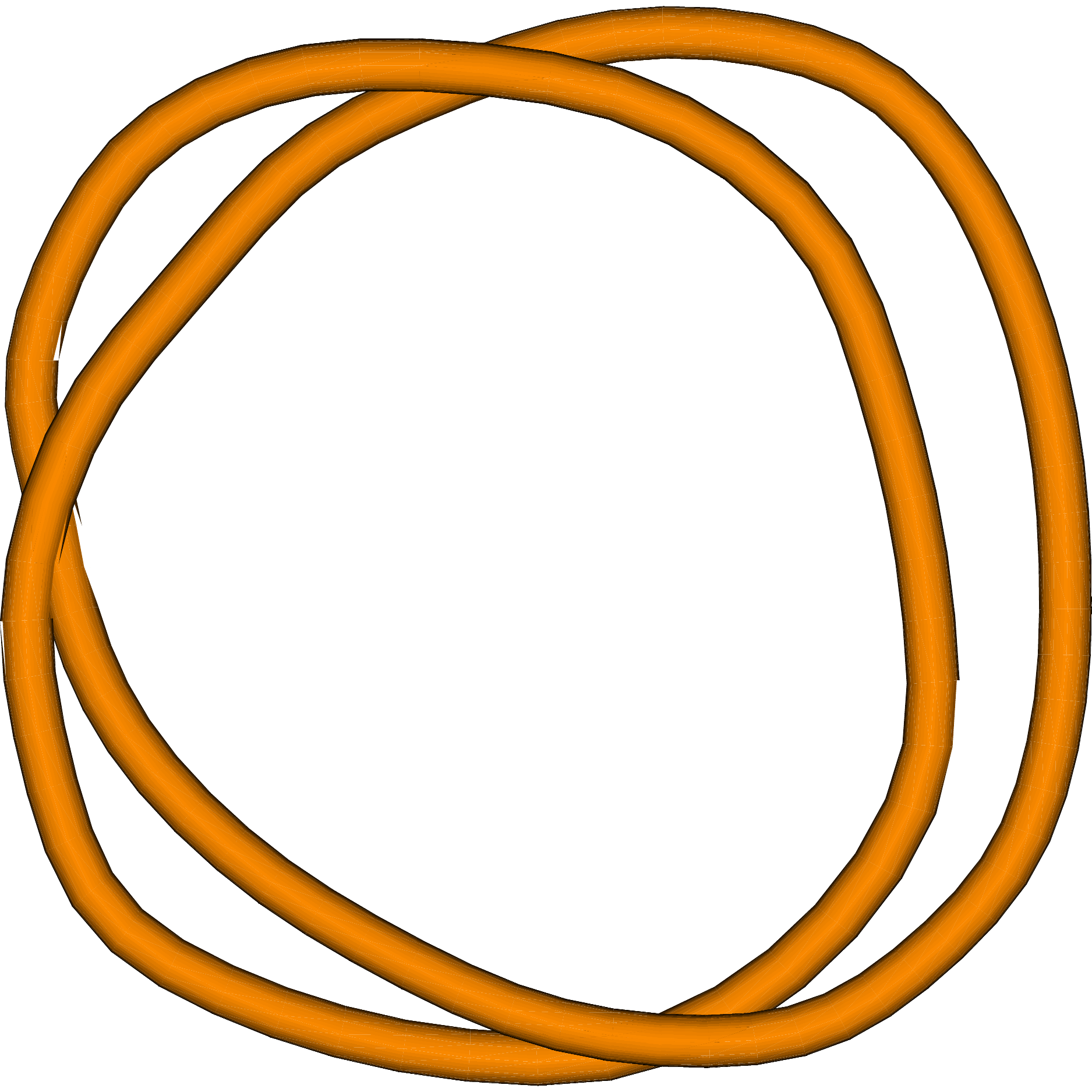}
\includegraphics[width=4.7cm]{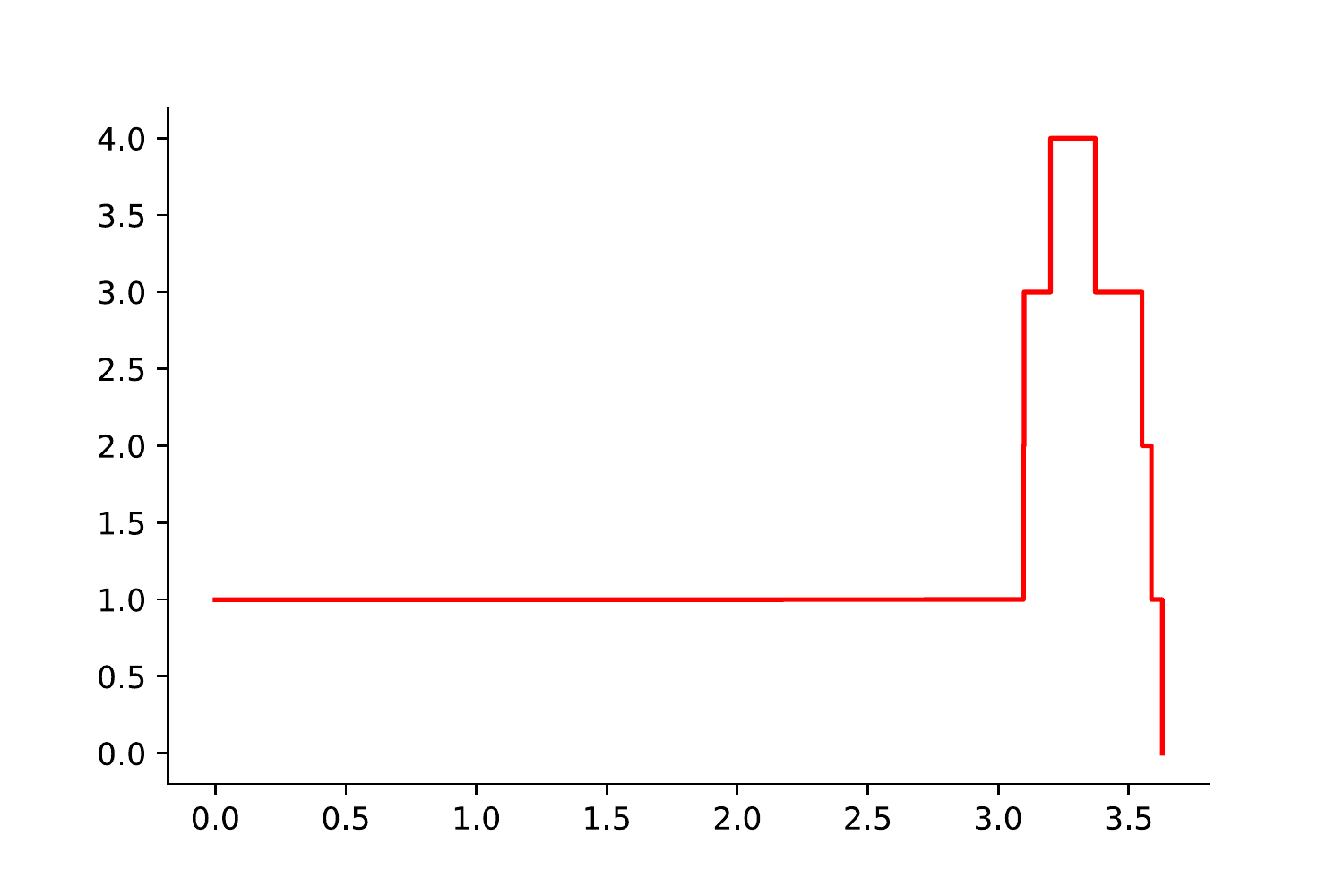}
\includegraphics[width=4.7cm]{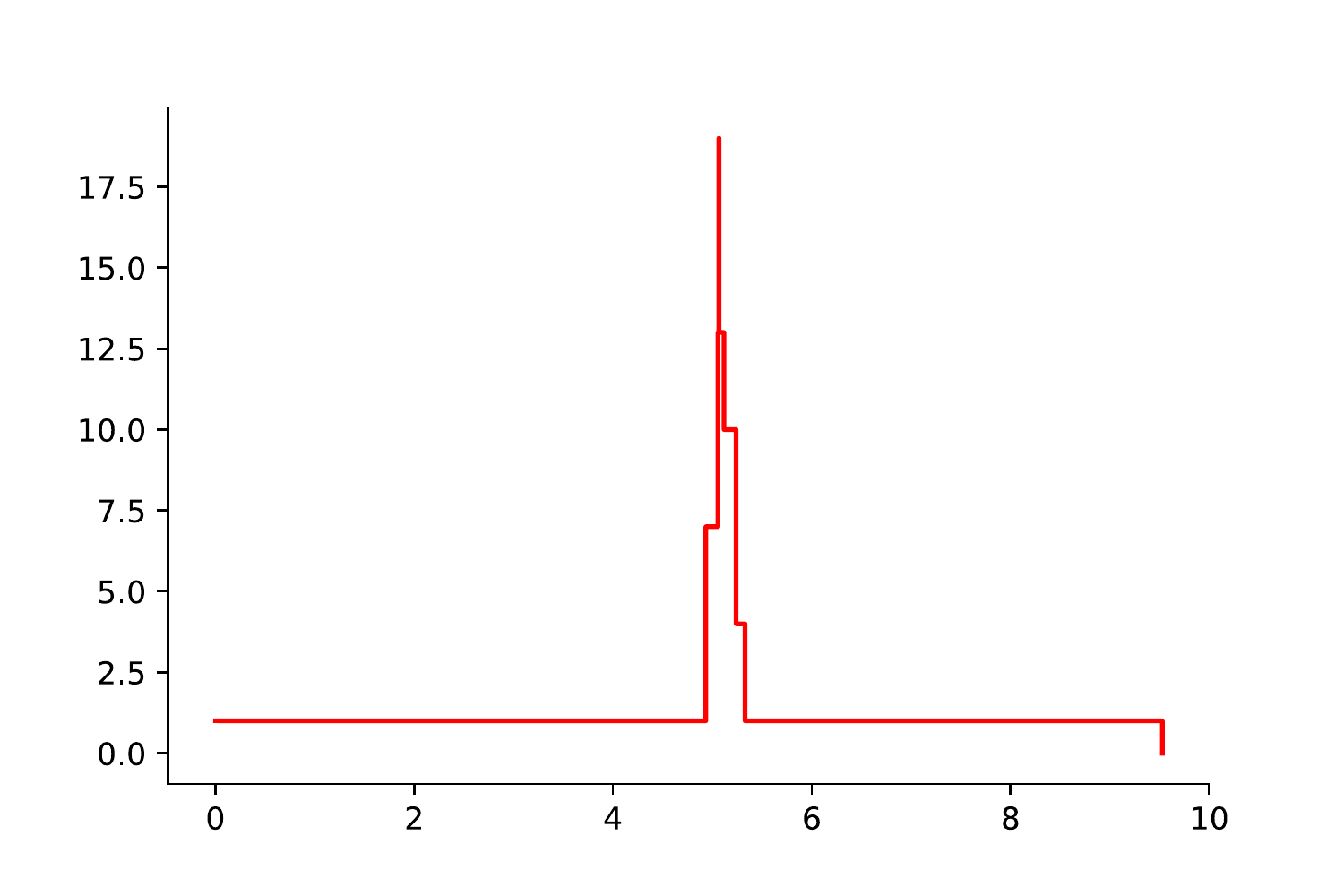}
\includegraphics[width=4.7cm]{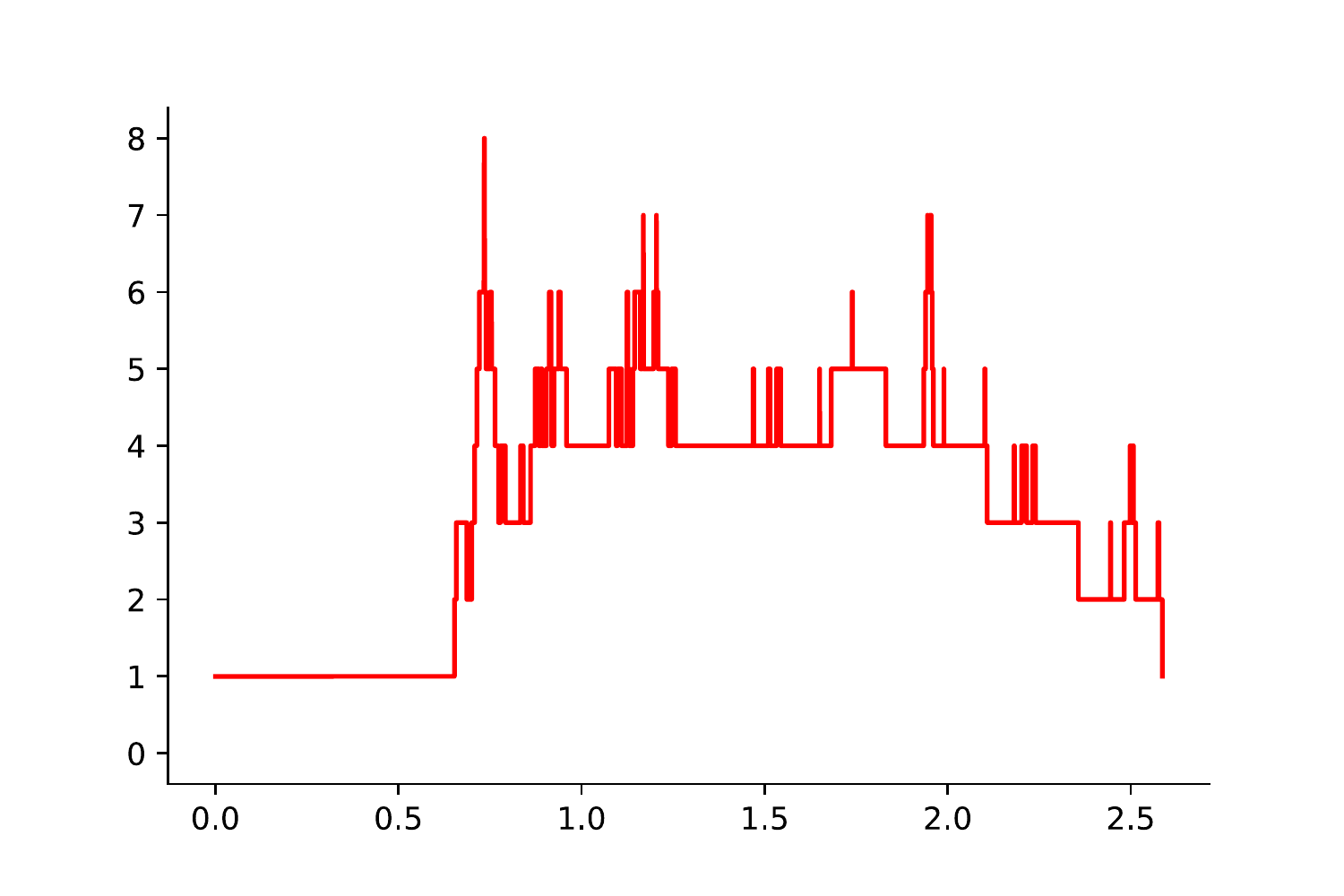}
\caption{Three different embeddings of the trefoil, and their corresponding (approximated) Betti curves. From left to right: a trefoil embedding close to its ideal configuration, a trefoil lying on a standard torus whose longitude is longer than its meridian, and a trefoil close to being planar. The images in the top panel were obtained with Knotplot~\cite{knotplot}.}
\label{fig:3trefoils}
\end{figure}
\end{exa}

\begin{rmk}~\label{rmk:ideal}
If $K$ is an ideal configuration~\cite{stasiak1998ideal} for the knot type $\mathcal{K}$, then we expect $\beta_1(K)(t)$ to be $1$ for $t \in [0,\frac{1}{2}[$ (since the injectivity radius is by definition $\frac{1}{2}$), and then jump to a large number $m(K)$ immediately after (see the left part of Figure~\ref{fig:3trefoils}). The number $m(K)$ is related to the number of self-tangencies of the ropelength minimiser embedding considered. We can use this to provide a measure of the ``closeness'' between a given embedding and  ideal configurations. We will define such a measure in Section~\ref{sec:distanceideal}.
\end{rmk}

\section{Persistent homology}\label{sec:PH}

\begin{figure}[ht]
\centering
\includegraphics[width=11cm]{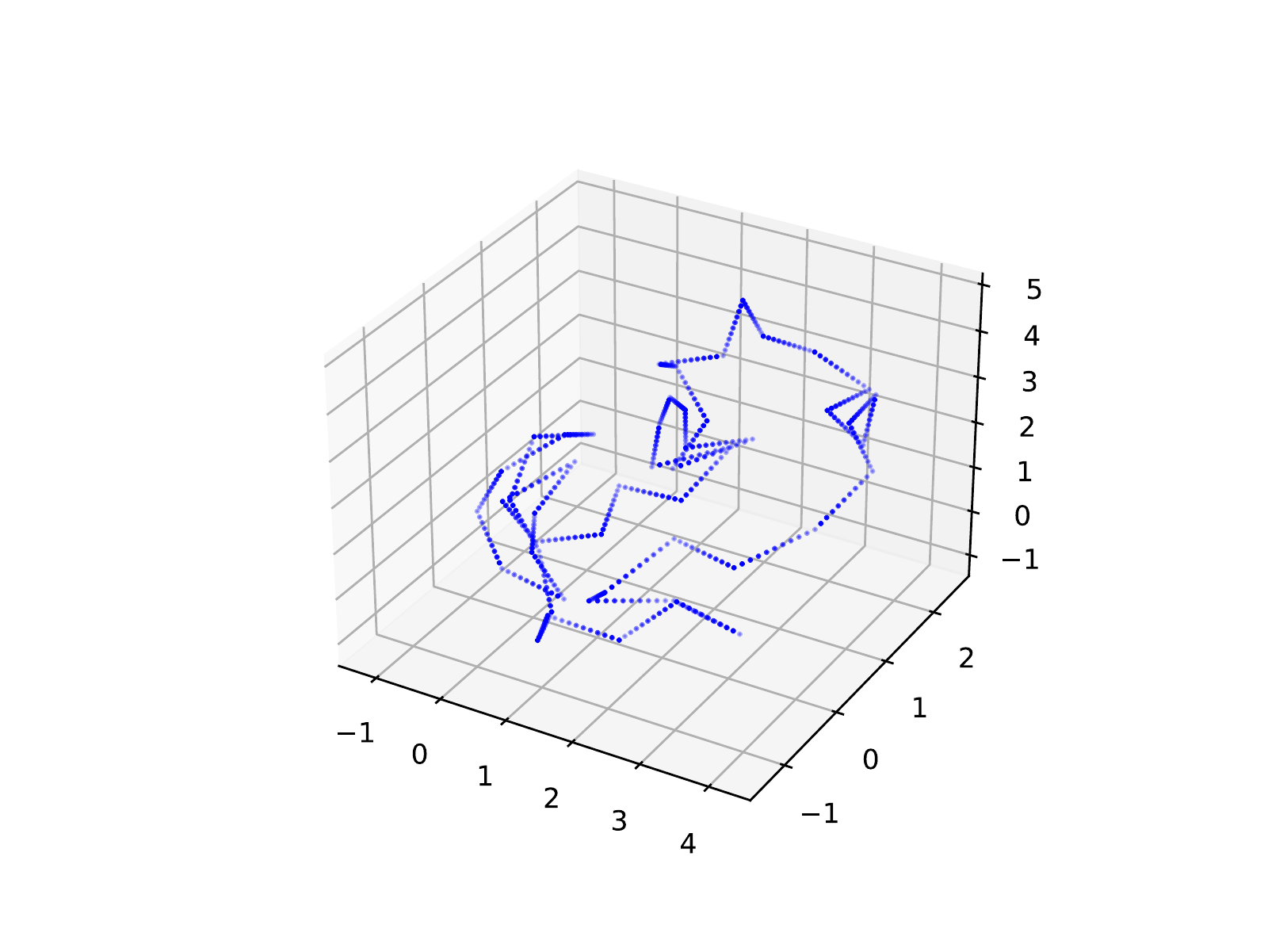}
\includegraphics[height=5cm]{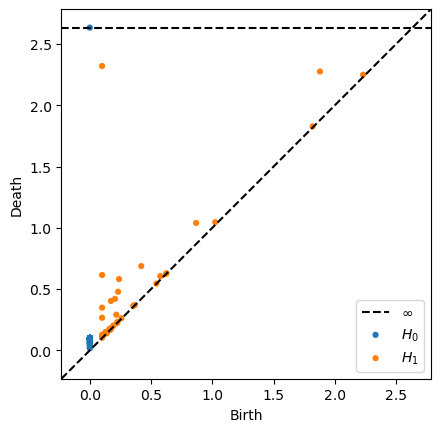}
\includegraphics[height=5cm]{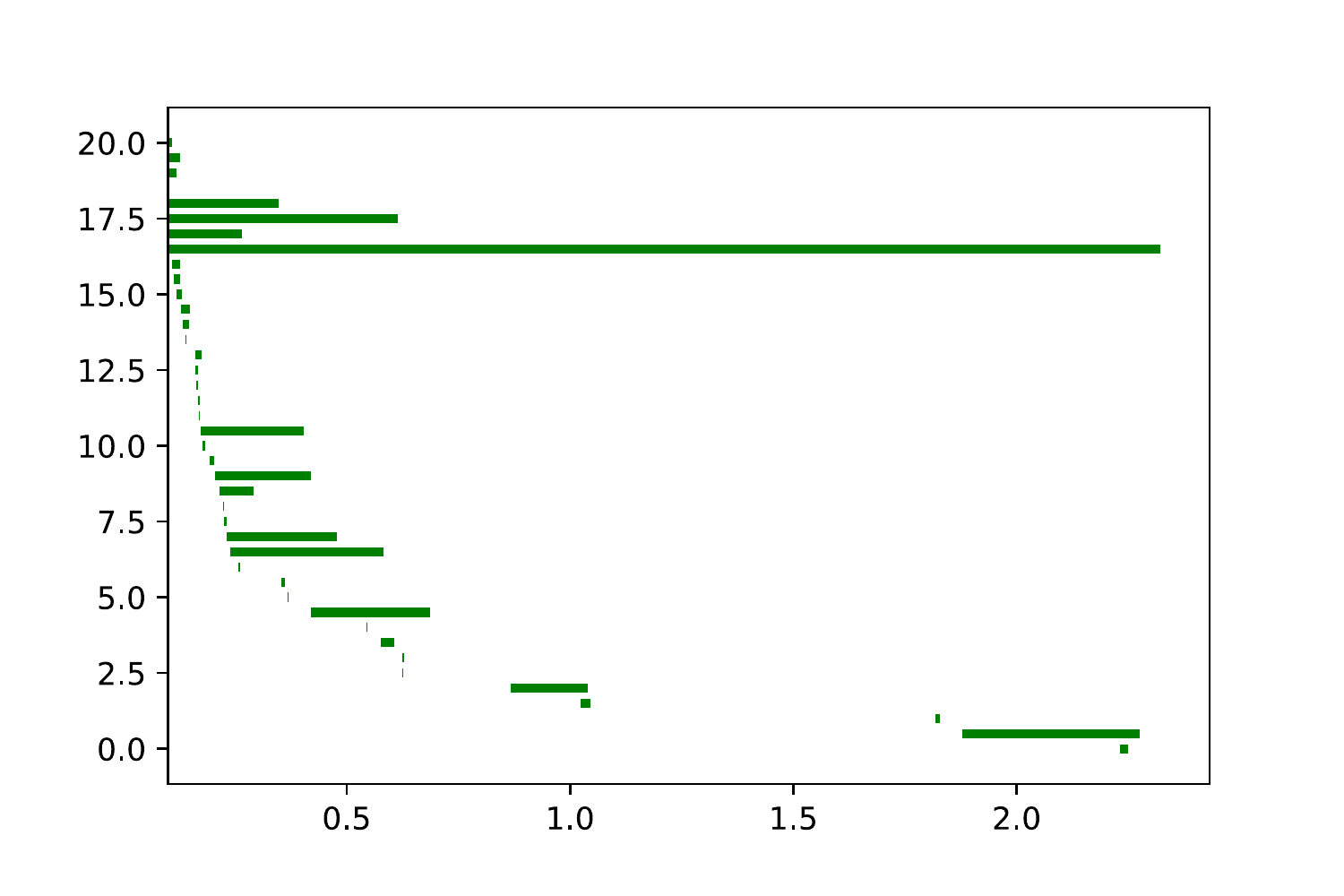}
\caption{In the top panel, a PL knot embedding of the trefoil with $50$ edges. Each edge is replaced with ten equidistant dots to create a point cloud (see also Figure~\ref{fig:representation}). In the bottom panel, two different ways to visualise the persistence diagram of this point cloud (in the barcode on the right we are only considering homology in dimension one). In the persistence diagram on the left, each point corresponds to a homology and its placement indicates the value of $t$ at which the class appears (birth) and the value of $t$ at which the class disappears (death). In the barcode on the right, homology classes are represented by bars that start and end at the birth and death value of $t$, respectively.}
\label{fig:PHdiagrams}
\end{figure}

Persistent homology \cite{edelsbrunner2010computational} is an algebraic tool for detecting topological and geometric properties of a space at different resolutions. 
The input for persistent homology is a nested sequence of simplicial complexes, called a \emph{filtration}. There are various ways to build a filtration on point-cloud data. Commonly used constructions are made from \v{C}ech complexes, Vietoris--Rips complexes, or $\alpha$-complexes. Roughly speaking, persistent homology captures the evolution of the homology of the filtered complex as it grows through the filtration. In particular, it keeps track of how homology classes appear and disappear, and this information can be represented in a \emph{persistence diagram} or \emph{barcode} (see Figure~\ref{fig:PHdiagrams}).  
The \v{C}ech complex built on a point cloud $P$ for a given radius $t \ge 0$ is the simplicial complex whose $0$-simplices are the points in $P$, and whose higher dimensional simplices are subsets of points in $P$ whose closed $t$-balls have non-empty intersection. This complex has the same homotopy type as the union of the closed $t$-balls centred in the points in $P$, by the \emph{Nerve Lemma} (see \emph{e.g.}~Section III.2 of \cite{edelsbrunner2010computational}). The \v{C}ech filtration on $P$ is the filtration consisting of the \v{C}ech complexes on $P$ for growing radii $t \geq 0$, and it gives a topologically faithful representation of the gradual thickening of the underlying space if $P$ is a sufficiently dense and uniform sample of the space.

However, for large data it can be computationally intensive --- and hence impractical --- to build a \v{C}ech filtration \cite{roadmap}. Therefore, in practice one may want to instead work with the Vietoris-Rips filtration, which depends only on pairwise distances between points and can therefore be computed much more efficiently. The Vietoris--Rips complex built on a point cloud $P$ for a given radius $t\geq0$ is the simplicial complex whose $0$-simplices are the point in $P$, whose $1$-simplices are the pairs of points in $P$ that are within a distance of $2t$ from each other, and whose higher dimensional simplices are the cliques of its $1$-simplices. The Vietoris--Rips filtration on $P$ is the sequence of Vietoris--Rips complexes for growing radii $t\geq 0$, and it is a good approximation to the \v{C}ech filtration in Euclidean space \cite{de2007coverage}.  

We use the Vietoris--Rips filtration on a point cloud $P(K)$ constructed from a PL knot embedding $K$ to approximate its metric neighbourhood for growing radii. 

We consider several features of the barcode corresponding to the Vietoris--Rips filtration on $P(K)$; the first are simply the number of bars and the length of the longest bar in the barcode in degree $1$ of the persistent homology of $P(K)$. These are denoted by $\# B(K)$ and $M(K)$ respectively. Furthermore, similarly to the Betti curve of a growing neighbourhood of a knot (see Definition~\ref{def:betticurve}), we can define a Betti curve on the rank of the homology of the Vietoris--Rips filtration on a point cloud. We compute the integrals of such Betti curves $\int_0^\infty \beta_1(P(K)) \mathrm{d}t$ by summing the lengths of all bars in the corresponding barcode. We denote the integral of the Betti curve corresponding to the Vietoris--Rips filtration on a point cloud sampled from a knot $K$ by $\mathcal{I}(K)$.  
We will see in Section~\ref{sec:results} that the integral $\mathcal{I}(K)$ as defined here constitutes a good approximation (with a few caveats, as explained in Figure~\ref{fig:intersectionangle}) of the integral of the Betti curve from Defintion~\ref{def:betticurve}.

\section{The data}\label{sec:data}

In this section, we describe the samples we study, as well as the methodologies used to generate and analyse them. 

All the random knots in our data sets are generated using the excellent Python based Topoly \cite{dabrowski2019topoly}. The same program is also used to determine the knot types of the PL curves considered whenever required. We use pyknotid~\cite{pyknotid} to compute average crossing numbers.
We use custom-made programs to compute the radius of gyration, curvature, torsion, and volume of the smallest enclosing sphere of a PL-embedded knot, and we use SciPy's built-in function to compute the volume of its convex hull.

We produce two qualitatively different data sets. For the first data set, we generate $10^4$ random knots for each length from $10$ to $100$ (in steps of $10$). We then compute the volume of the minimal enclosing sphere and that of the convex hull inscribing each knot, as well as the curvature, torsion, ACN and radius of gyration. 

For the second data set, we sample $10^3$ knots for each knot type with up to $6$ crossings, for lengths between $50$ and $200$ (in increments of $50$); we also include $10^3$ samples of knots whose knot types do not belong to the previous categories (so whose minimal crossing number is $\ge 7$), and we will refer to these as ``unknown'' in what follows.

\begin{figure}[ht]
\centering
\includegraphics[width =6.5cm]{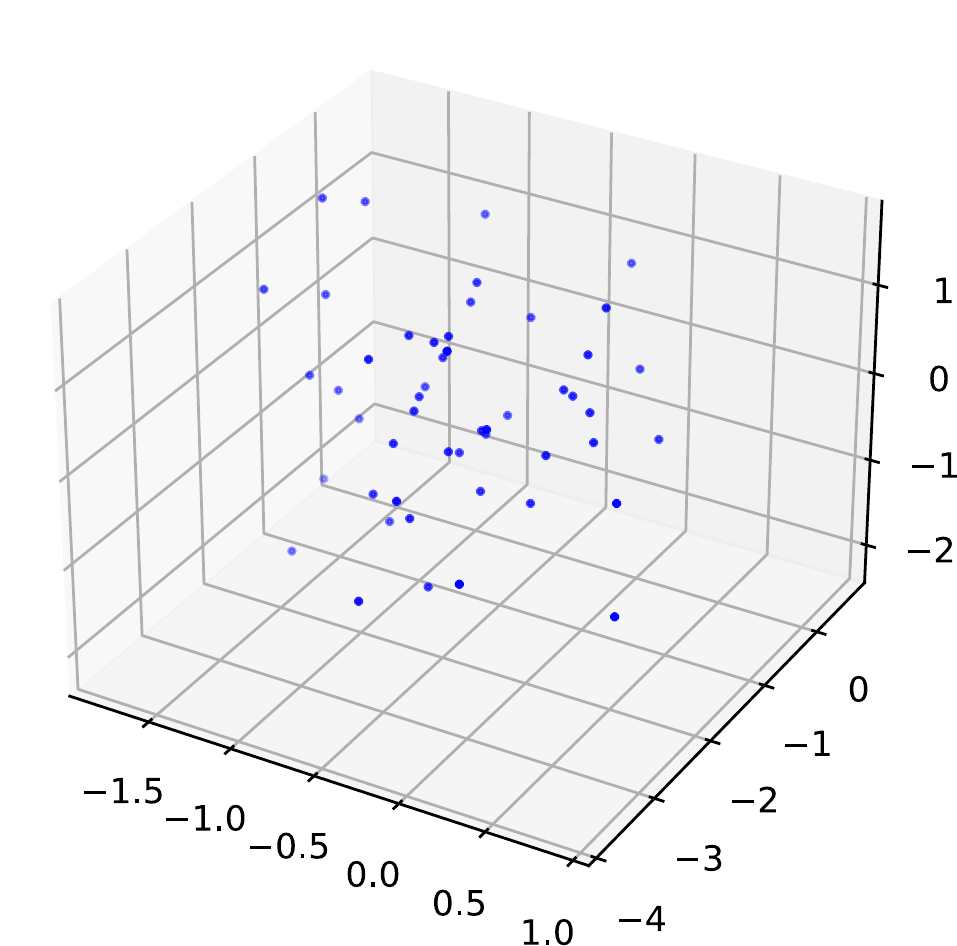}
\includegraphics[width =6.5cm]{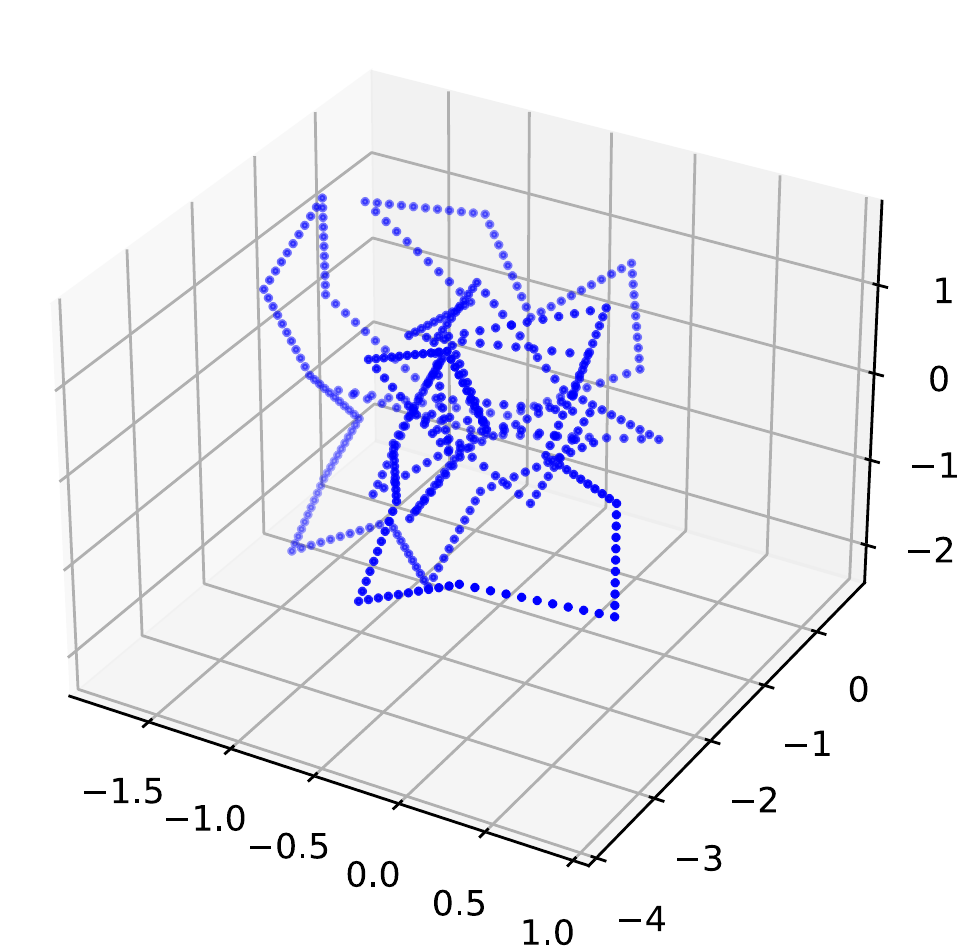}
\caption{A point cloud describing the endpoints of a PL embedding $K$ of the trefoil with $50$ edges, and the point cloud $P(K)$ we consider, obtained by interpolating the unit length edges with ten further equidistanced points.}
\label{fig:representation}
\end{figure}

We then compute the barcodes for Vietoris--Rips filtrations associated to the knots as explained below using the efficient program Ripser~\cite{ripser}. As our purpose is to approximate the topology of the neighbourhoods of our generated embeddings as closely as possible, we do not simply compute the PH of the Vietoris--Rips filtrations on the endpoints of the sampled PL curves.
Instead, we interpolate the endpoints of the unit segments of each embedding with ten equidistant points, and compute the Vietoris--Rips filtration on the resulting point cloud, which we denote by $P(K)$ (see Figure~\ref{fig:representation}). Using ten points per linear segment gives a dense enough sample of the embedding to yield a Vietoris--Rips filtration that approximates the growing tubular neighbourhood closely for most knots.
For practical reasons, we restrict our considerations to the first homology groups. It is likely that higher homology groups do also retain useful information on the geometric structure of such embeddings.

In order to be able to find meaningful correlations between geometric aspects --- such as enveloping volumes, curvature, torsion, and the ACN --- and PH, we consider the features extracted from the barcodes described in the previous section.
Note that, in the case at hand, the support of $\beta_1(K)$ is contained in the interval $[0, RS(K)]$, where again $RS(K)$ is the diameter of the minimal  sphere enclosing the embedding. Indeed, for $t\ge RS(K)$ the neighbourhood $\nu_t(K)$ is topologically a $3$-ball. This readily implies that $\mathcal{I}(K)$ is well defined, \emph{i.e.}~it is a proper integral.

The correlations displayed in the next section were computed using SciPy's statistics module. All the programs we used to generate the data, as well as the data itself, is available on the first author's GitHub page~\cite{miogithub}.

\section{Results}\label{sec:results}

In this section, we collect the results of the computations detailed above.
More precisely, we show that the correlation between the integral of $\beta_1$  and the volume occupied by a knot becomes increasingly negative as the knot's length increases. The magnitude of these correlations is especially large in the case of the volume of the convex hull.
This confirms the intuitive fact that, for long embeddings, being geometrically complex implies being spatially compact. Furthermore, we quantify how considering increasingly complex topologies (\emph{i.e.}~progressively complex knot types) influences this inverse correlation. 
\begin{figure}[ht]
\centering
\includegraphics[width=6.5cm]{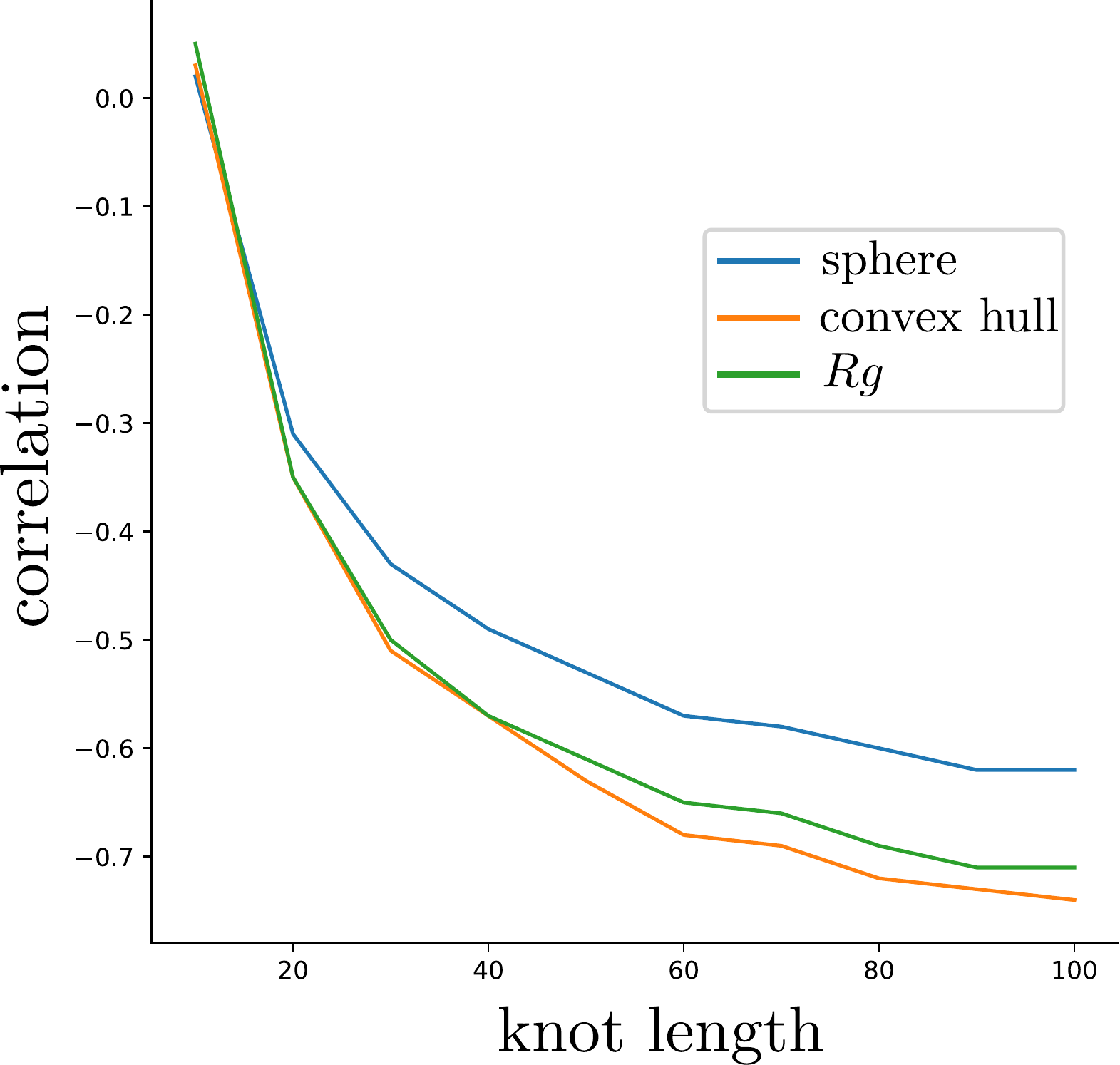}
\includegraphics[width=6.5cm]{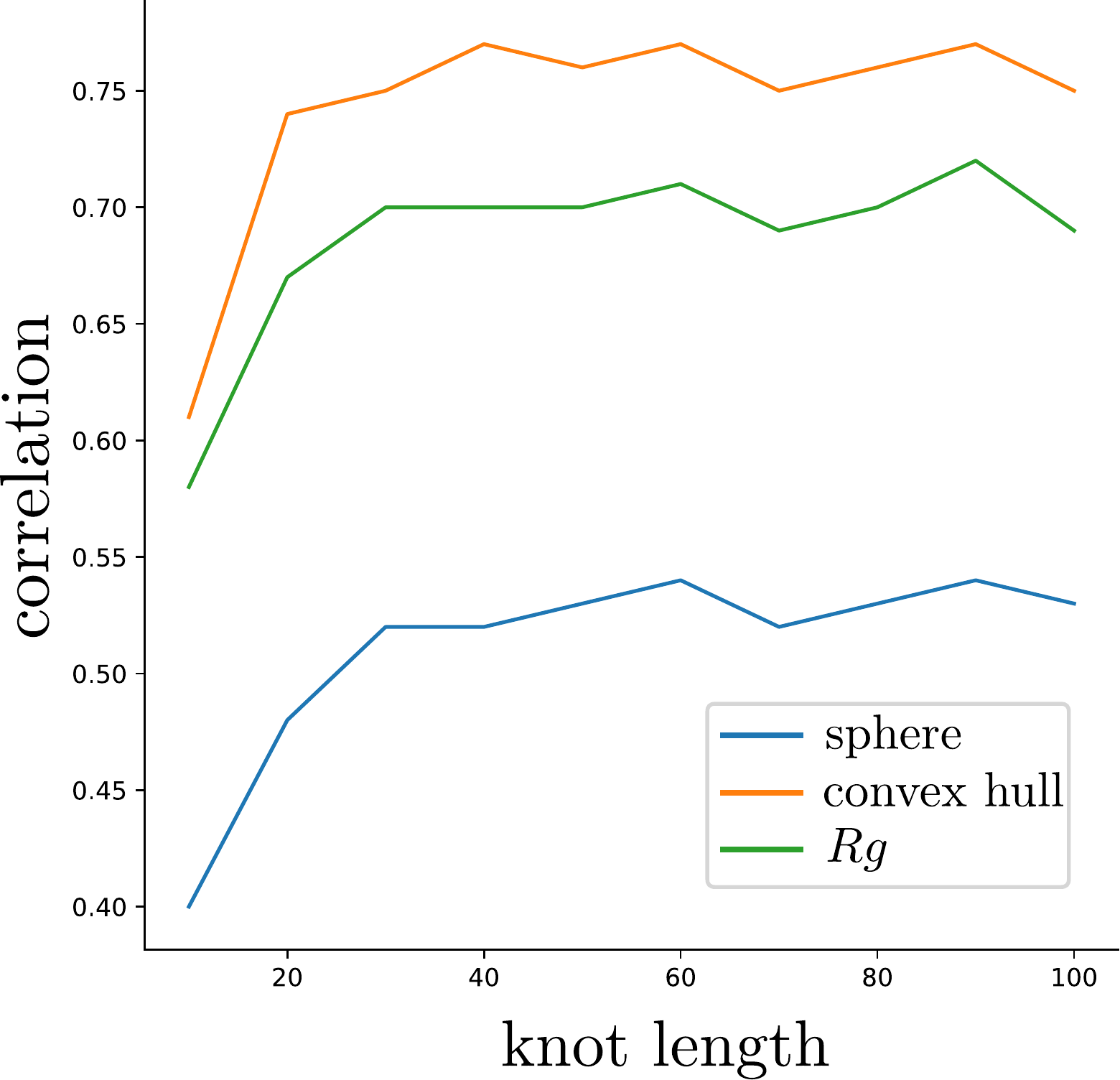}
\caption{Left: the correlation between $\mathcal{I}(K)$ and the volume of the circumscribing sphere, convex hull and $Rg$ as a function of  knot length. Right: the correlation between the length of the maximal bar $M(K)$ in the Vietoris-Rips barcode and the same quantities as above.}
\label{fig:correlsINTRIPS}
\end{figure}

\begin{figure}[ht]
\centering
\includegraphics[width=6.5cm]{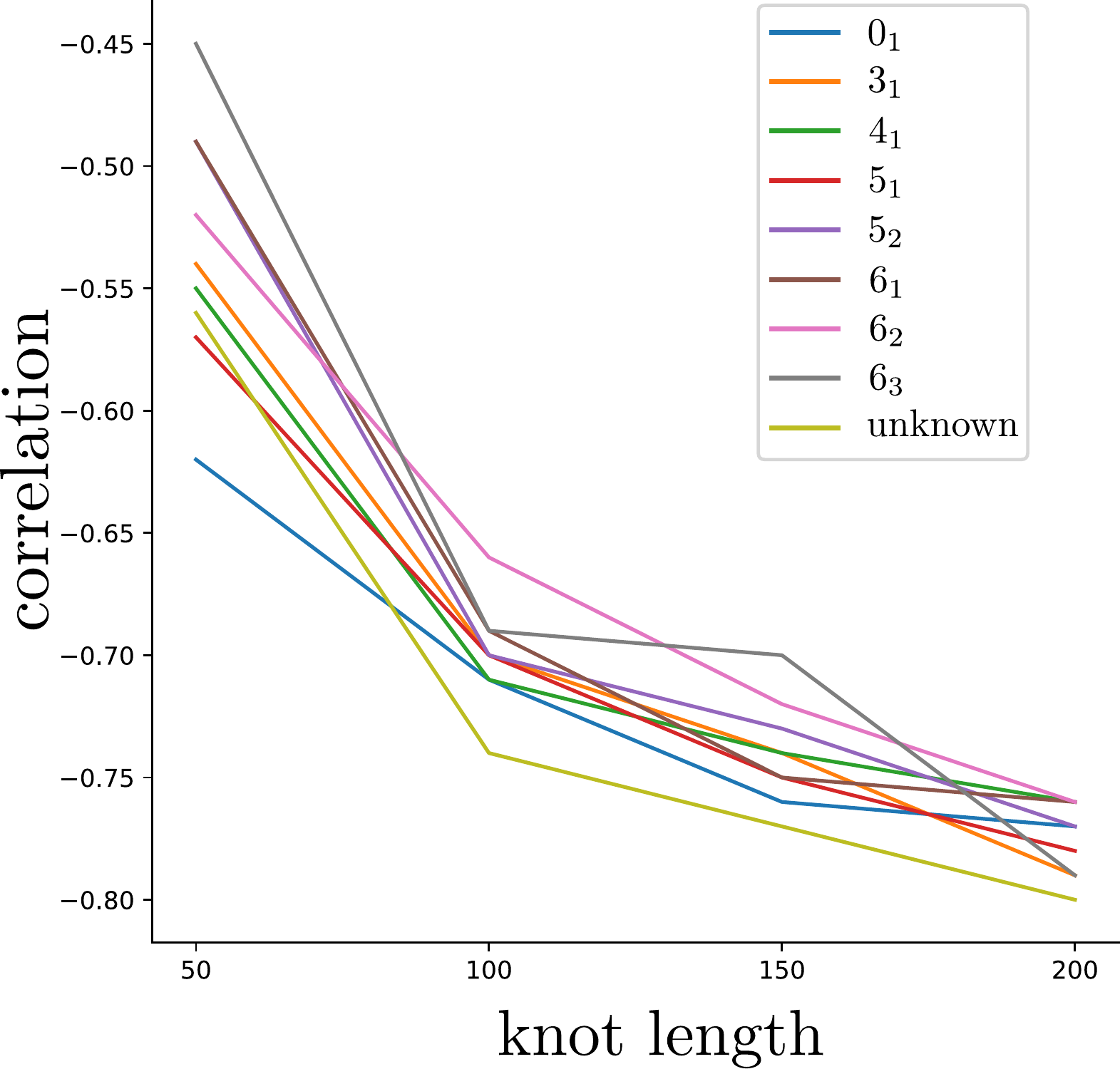}
\includegraphics[width=6.5cm]{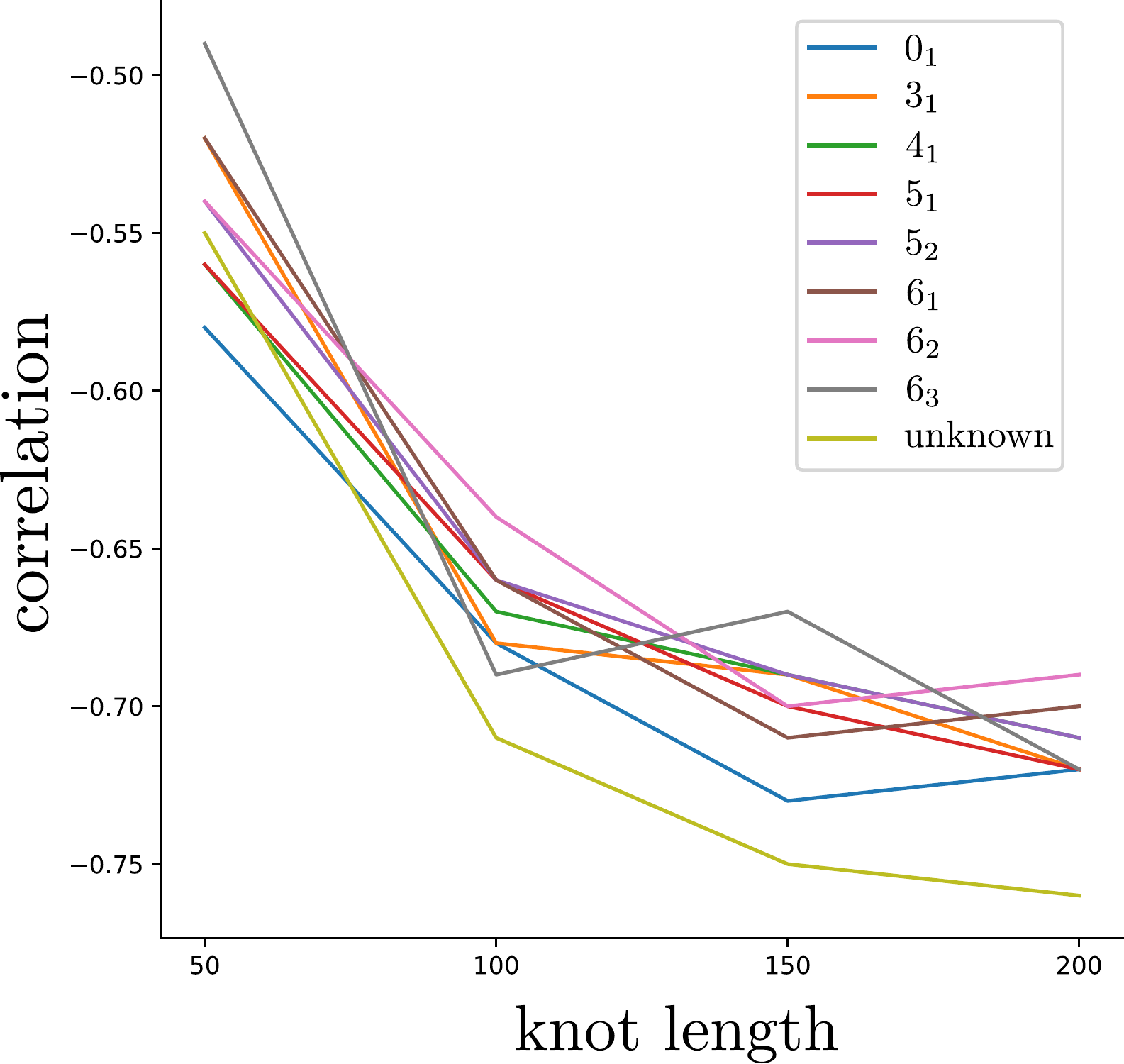}
\caption{Here we focus on the correlations between $\mathcal{I}(K)$ and the volume of the convex hull (left) and $Rg$ (right).  We are splitting the correlations according to the various knot types considered.}
\label{fig:ALLcorrelsINTRIPS}
\end{figure}

\begin{figure}[ht]
\centering
\includegraphics[width=6.5cm]{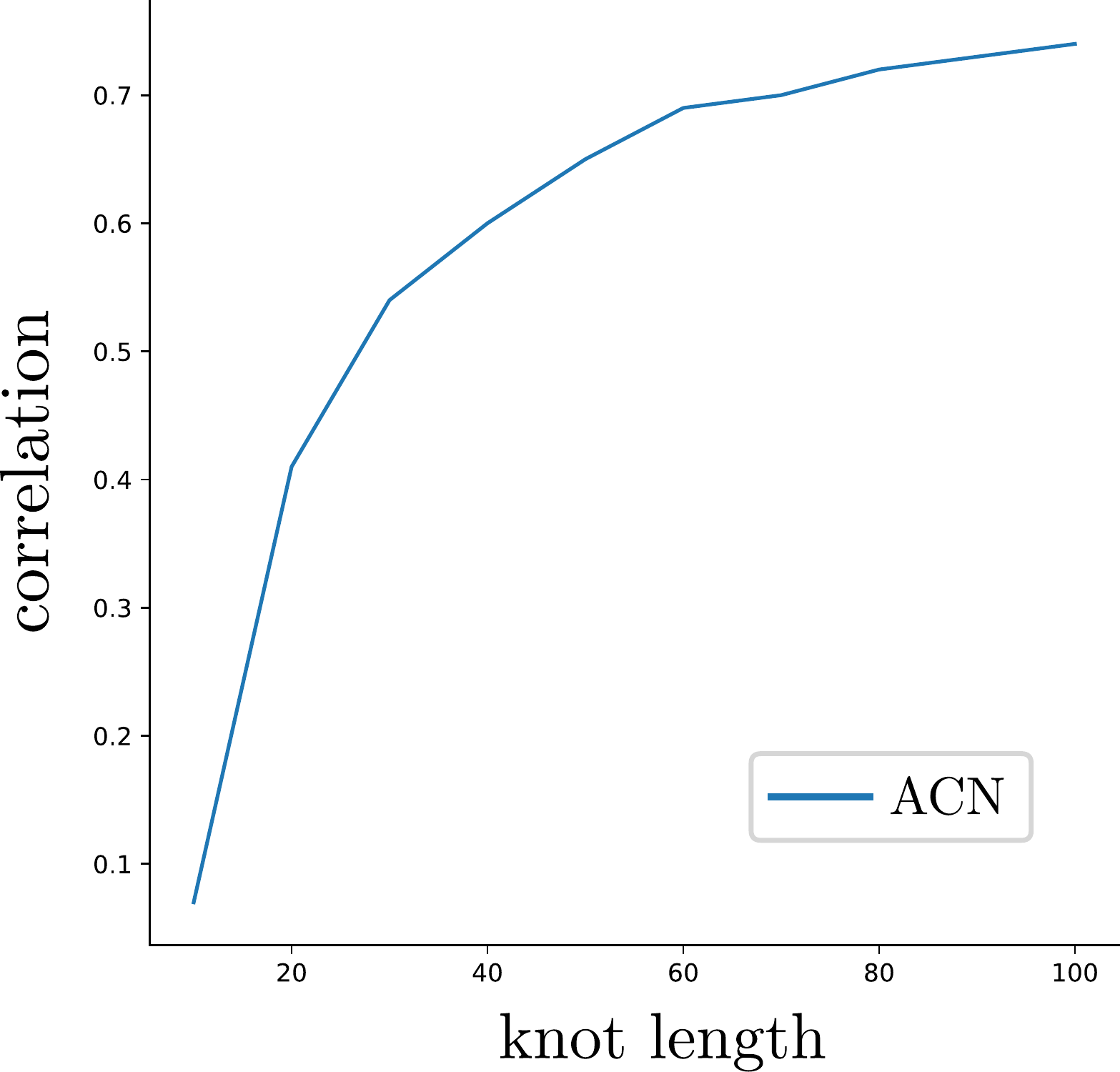}
\includegraphics[width=6.5cm]{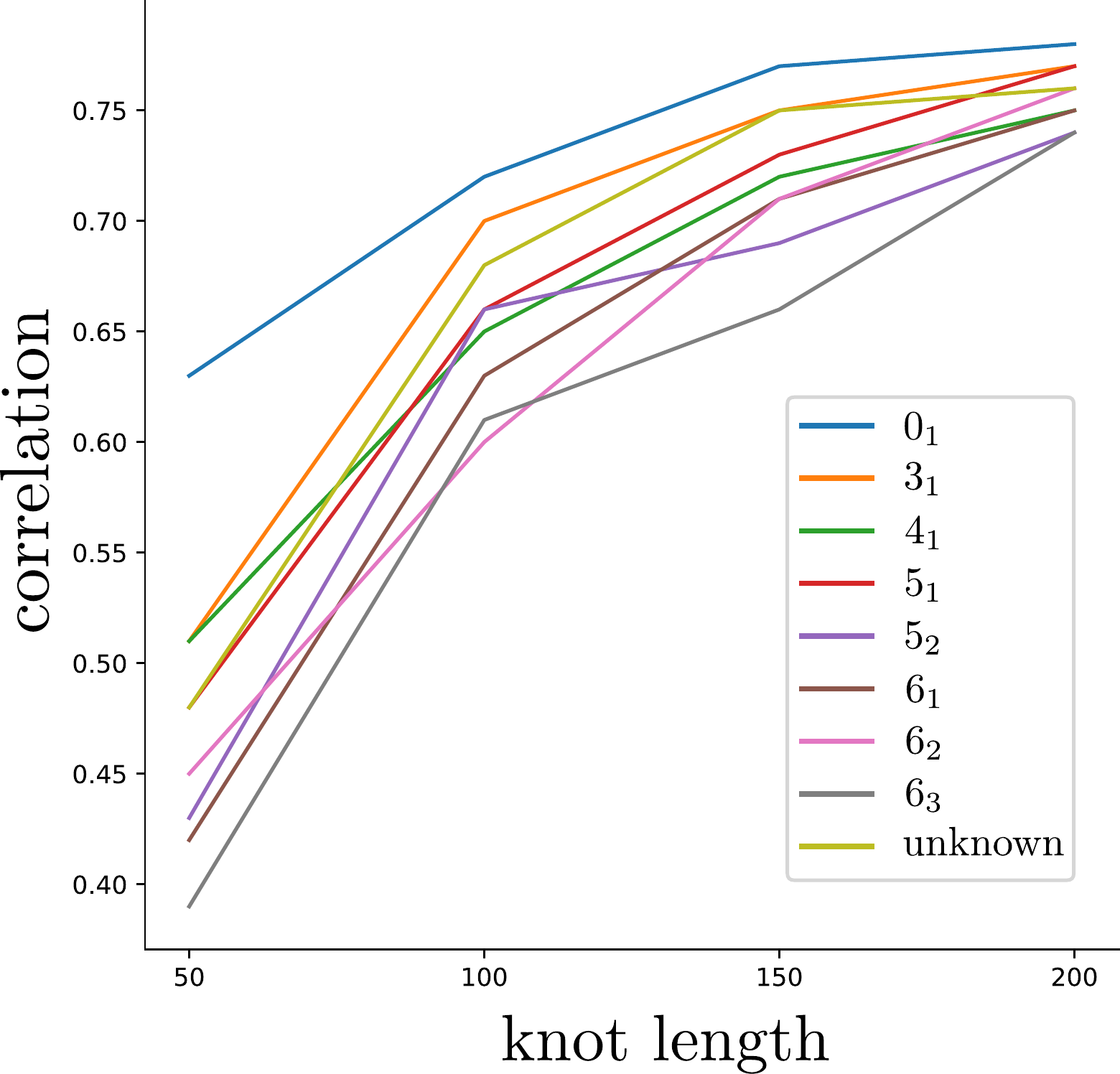}
\caption{Correlation functions between the average crossing number and $\mathcal{I}(K)$. On the left we consider all knots (with lengths ranging from $10$ to $100$), while on the right we are splitting among the various knot types considered (with lengths ranging from $50$ to $200$).}
\label{fig:corrsACN}
\end{figure}

Similarly, we find a growing correlation between $\mathcal{I}(K)$ and the average crossing number (left panel of Figure~\ref{fig:corrsACN}). This is not surprising, as both can be thought of as being measures of the geometric complexity of the embedding; it is however interesting to observe the phenomena displayed in the right panel of Figure~\ref{fig:corrsACN}, where this correlation's behaviour is split among the different knot types. Here we observe that, rather unexpectedly, the values of the plots do not appear to be monotonically related to the complexity measure on knot types given by the minimal crossing number. Indeed, the ``unknown'' category has values which are larger than most other knot types (in the range of lengths considered). At the same time, the unknots' correlations appear to be  considerably larger. 

We also report the existence of an increasingly positive correlation between $M(K)$ (as displayed the right panel of Figure~\ref{fig:correlsINTRIPS}) and the volume of the convex hull. Other correlations between~\emph{e.g.}~$\#B(K)$ and the convex hull's volumes or $Rg$ are marginally weaker than what we observe for $\mathcal{I}(K)$, and we did not include this data in this manuscript, referring to the full data set in~\cite{miogithub}. \newline

Interestingly, unlike in the case of the average crossing number, we do not find any significant correlation between either the curvature or torsion of the knots and the integral $\mathcal{I}$. This should be compared with other more ad hoc approaches, such as~\cite{curvtors,acnandwrithe,totcurvtors}.\newline

\begin{figure}[ht]
\centering
\includegraphics[width=12cm]{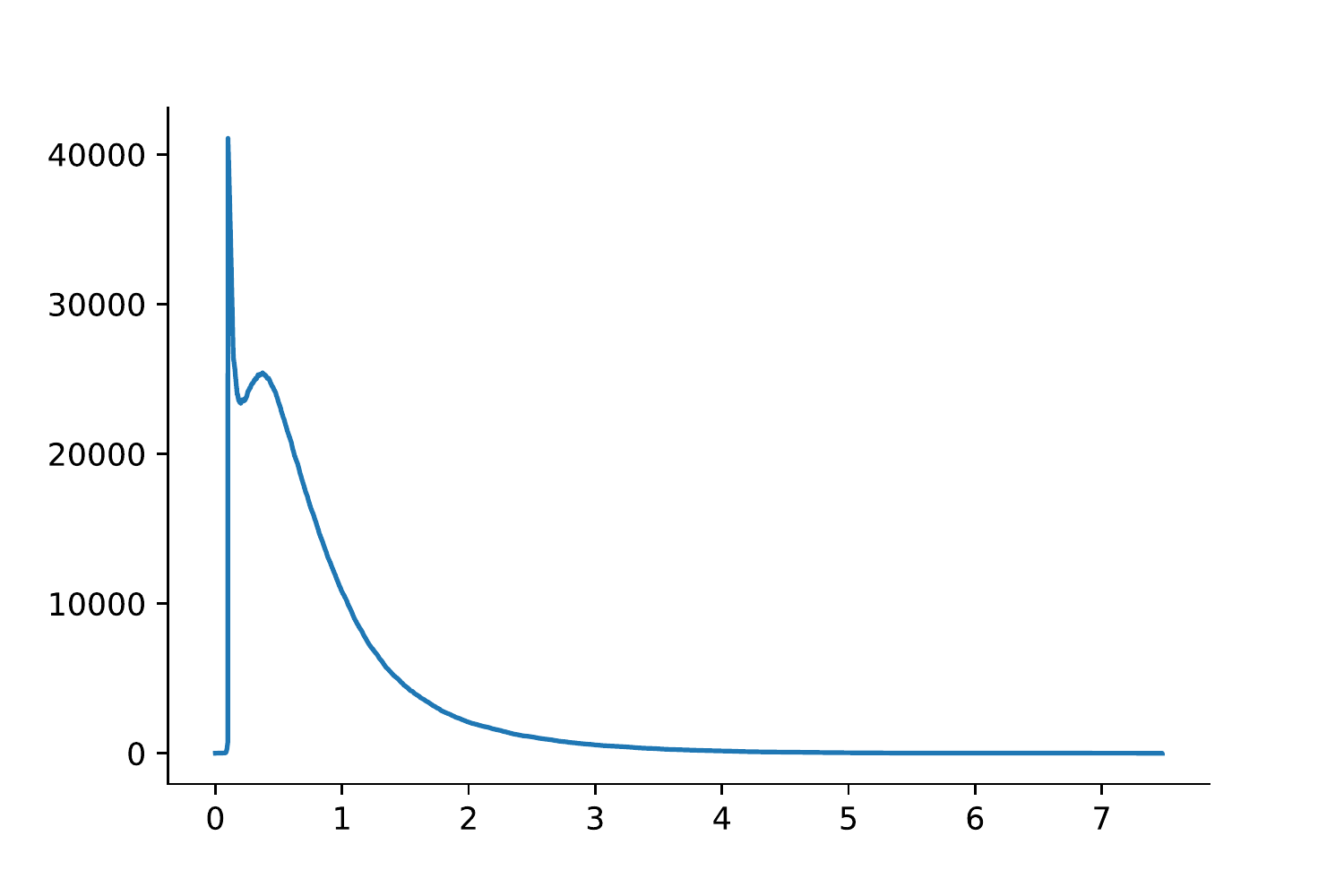}
\caption{The sum of $1000$ Betti curves for random unknots of length $200$. This overall shape is present for all knot types considered. The average curve can be obtained simply by dividing by $1000$ pointwise.}
\label{fig:average-0_1}
\end{figure}

We also look at the sum of all Betti curves for each length and\slash or knot type. It turns out that the overall shape of this ``cumulative'' curve is the same, regardless of the knot type.  One example of such a curve is shown in Figure~\ref{fig:average-0_1} for the case of  unknots of length $200$. Note that the parameter $t$ used by Ripser is the diameter of the points' neighbourhoods, rather than the radius, which has the effect of stretching the domain's length by two in Figure~\ref{fig:average-0_1}.

\begin{figure}[ht]
\centering
\includegraphics[width=8cm]{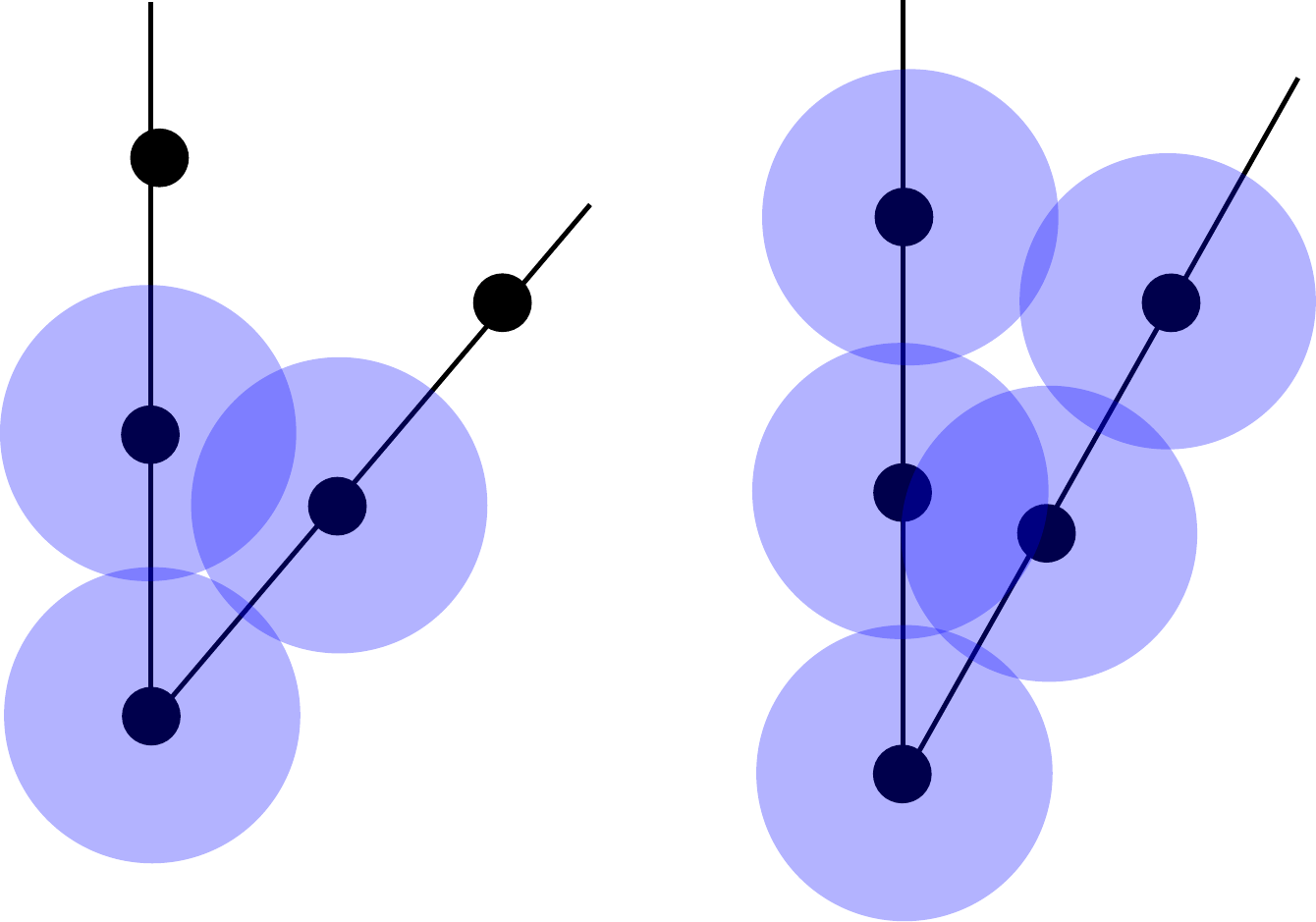}
\caption{Left: A short lived homology class in the \v{C}ech complex of $P(K)$; these appear for angles smaller than $\frac{\pi}{3}$, when the vertex point of a PL knot embedding $K$, together with the two adjacent points in $P(K)$, have neighbourhoods that only intersect in pairs for a small range of radii. When all three disks intersect, the length-three $1$-cycle they generate in the \v{C}ech complex is capped off by a $2$-simplex. Right: the analogous situation for the Vietoris-Rips complex, where four vertices are needed to create a short lived class. In this case the angle has to be smaller than $\arccos{\left(\frac{3}{4}\right)}$.}
\label{fig:intersectionangle}
\end{figure}

It is interesting to point out certain clear characteristics of these cumulative functions, which are present for each choice of knot type and length. The leftmost spike indicates the presence of a large amount of rather short lived homology classes that appears right after $t = 0.1$. This is due to the fact that two consecutive segments in a PL knot embedding whose internal angle is less than or equal to $\arccos{\left(\frac{3}{4}\right)} \sim 41.4^\circ$ will produce such a class, as shown in Figure~\ref{fig:intersectionangle}.
Therefore, the presence of this spike can be regarded as a consequence of our choice to use the Vietoris-Rips complex as an approximation of a knot's neighbourhood. However, as pointed out in the left part of Figure~\ref{fig:intersectionangle}, a similar (potentially wider and higher) spike would have appeared if we had chosen the \v{C}ech complex instead. 
Incidentally, the value of the maximum attained is therefore related to the average number of small angles present in the embeddings. It is possible to remove this spike by simply ignoring all short lived bars appearing right after $t =0.1$.

Similarly, the value of the second maximum might be of interest, as it appears to be related to the average distance between the edges of PL embeddings in the given knot type.

\begin{figure}[ht]
\centering
\includegraphics[width=6.5cm]{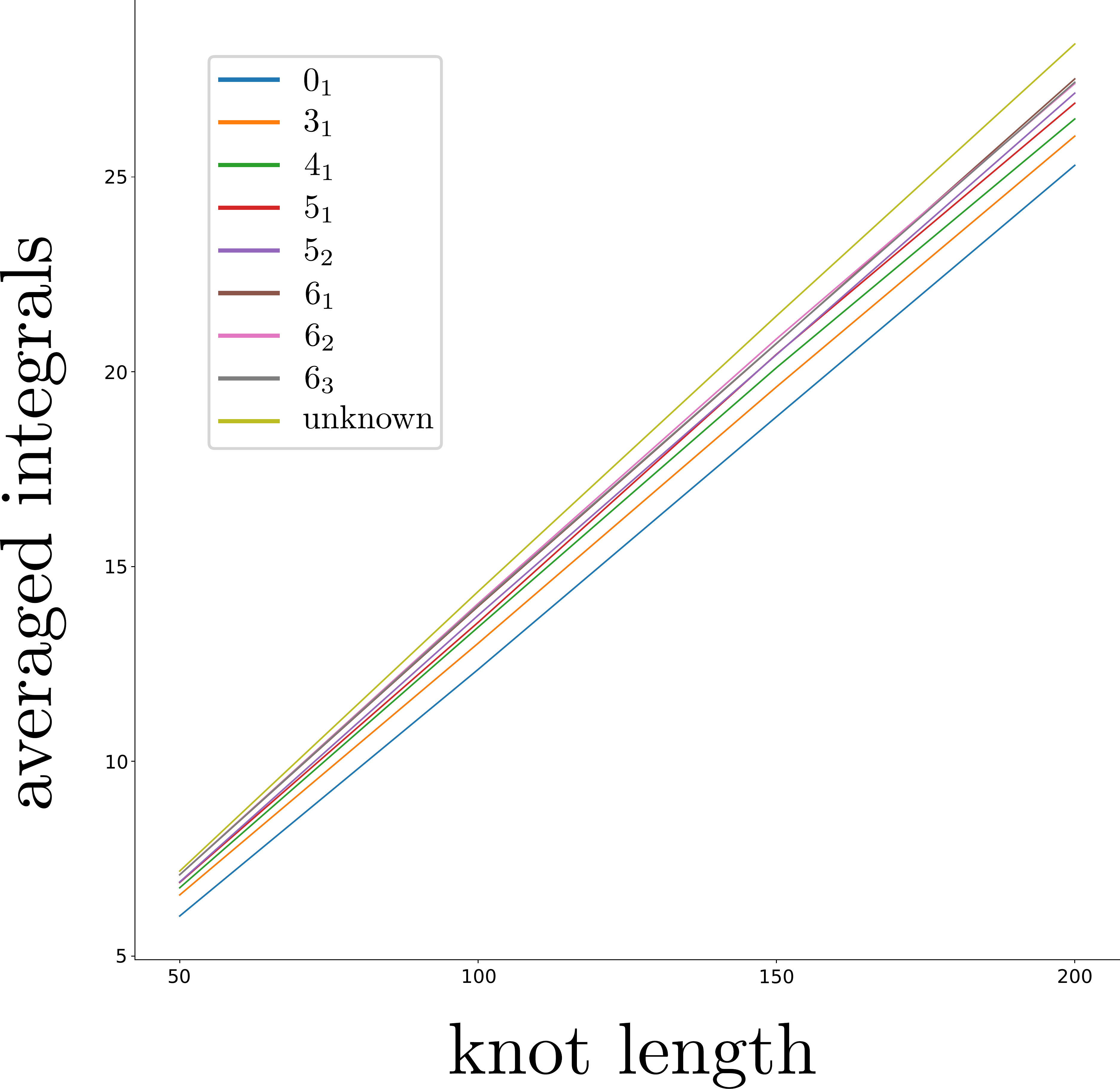}
\includegraphics[width=6.5cm]{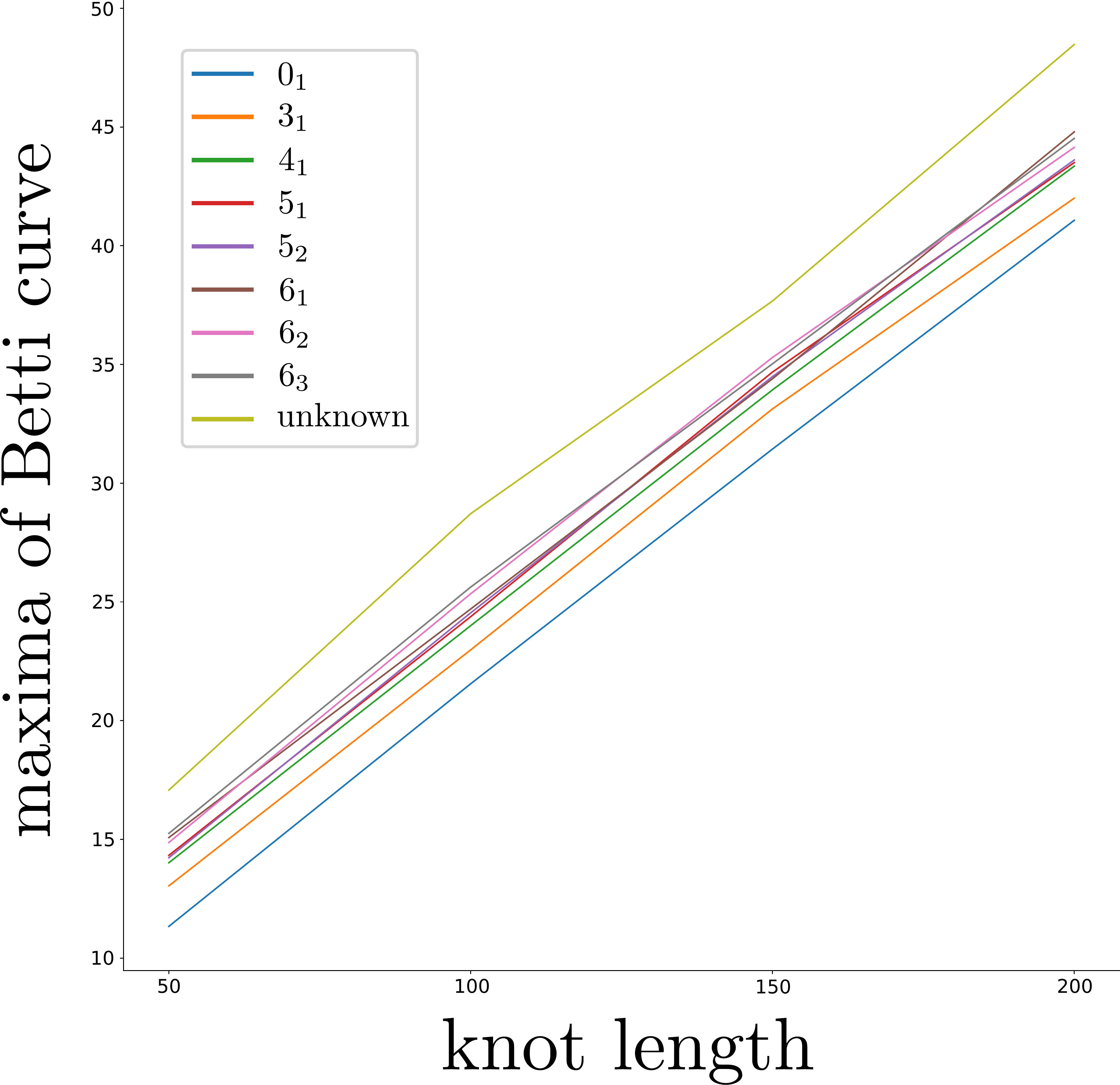}
\caption{Left: average of the integrals of the Betti curves for each knot type as a function of length.  Right: average of the maxima of the Betti curves (see Figure~\ref{fig:average-0_1}) for each knot type as a function of length.}
\label{fig:valuesintmaxima}
\end{figure}

We display the values of the average integral and the maximum of the average Betti curves for the various knot types considered as a function of length in Figure~\ref{fig:valuesintmaxima}. In both cases, there appears to be an almost perfect linear relation. Further, unlike what can be observed in Figures~\ref{fig:ALLcorrelsINTRIPS} and \ref{fig:corrsACN}, here we have a clear division into knot types, with values increasing monotonically with the knot type's topological complexity.

\section{Deviation from ideality}\label{sec:distanceideal}
In this final section we use the tools considered so far to introduce a naive way of quantifying ``how much'' a given knot deviates from being ideal. The key idea follows from Remark~\ref{rmk:ideal}: we can use the fact that the PH of ideal knots has a predictable behaviour to quantify the dissimilarity of the given knot from ideal knot embeddings. Note that our measure will quantify the difference of a given embedding to all possible ideal embeddings, rather than just those belonging to the same knot type.\newline

Denote by $S(K) = \max\{t \ge 0\,|\,\beta_1(K) \neq 0\}$, and let  $\varepsilon$ be a small positive real number. Call $f_{R,\varepsilon}: [0,R] \longrightarrow \mathbb{R}_{\ge 0}$ the unique function obtained by considering the linear function taking value $1$ on $0$, and value $0$ on $R-\varepsilon$, and defined to be identically $0$ after $R-\varepsilon$. Then define 
\begin{equation}\label{eqn:deltameasure}
\delta_\varepsilon (K) = \frac{1}{S(K)} \int_0^{S(K)} f_{S(K),\varepsilon}(t)\cdot\max\{\beta_1(K)(t) -1, 0\} \mathrm{d}t.
\end{equation}

We claim that $\delta_\varepsilon$ defines a sensible quantification of the ``distance'' between the given embedding and ideal ones. Let us examine the various components of Equation~\eqref{eqn:deltameasure} to substantiate the claim. Recall from Remark~\ref{rmk:ideal} that, for an ideal embedding $K_I$ of a knot, the function $\beta_1(K_I)$ takes the value $1$ until $t=\frac{1}{2}$, where it jumps to $m(K)$, indicating the appearance of $m(K)$ bars. These may not be the last bars appearing, but it is reasonable to assume that any further bar will be short-lived. This is because, by definition, the self-touching solid torus $\nu_{IR(K)}(K)$ ``occupies'' most of the volume surrounding the knot.

As we are modelling a knot, the value of $\beta_1(K)$ will usually be at least $1$ on the interval $[0,S(K)]$; we thus calibrate for this information by considering $\max\{\beta_1(K)(t)-1,0\}$. The role of the function $f_{R,\varepsilon}$ is to preserve the contribution of bars appearing early (\emph{e.g.}~for small values of the diameter), and to erase the contribution of bars appearing towards the end of the support of $\beta_1(K)$. Of course, this is not the only possible choice of a function with this property, but it is definitely one of the simplest.
The value of $\varepsilon$ acts as a cutoff, meaning that all bars born after $R-\varepsilon$ will not contribute to $\delta_\varepsilon$. 

It can be checked that as expected, for sufficiently small choices of the $\varepsilon$ threshold, for the trefoils in Figure~\ref{fig:3trefoils} we have (from left to right) $0<\delta_\varepsilon(T_1) < \delta_\varepsilon(T_2) < \delta_\varepsilon(T_3)$.


\begin{thebibliography}{1}
\bibliographystyle{plain}

\bibitem{packing}
Annunziato, A.
\newblock DNA packaging: nucleosomes and chromatin
\newblock {\em Nature Education} 1, no. 1 (2008): 26.

\bibitem{phages}
Arsuaga, Javier, Vázquez Mariel, Trigueros Sonia, and Roca Joaquim.
\newblock Knotting probability of DNA molecules confined in restricted volumes: DNA knotting in phage capsids
\newblock {\em Proceedings of the National Academy of Sciences } 99, no. 8 (2002): 5373-5377.


\bibitem{landscapes}
Bubenik, Peter.
\newblock Statistical topological data analysis using persistence landscapes
\newblock {\em J. Mach. Learn. Res.} 16, no. 1 (2015): 77-102.

\bibitem{carlsson2009}
Carlsson, Gunnar.
\newblock Topology and data
\newblock {\em Bulletin of the American Mathematical Society} 46, no. 1 (2009): 255--308.

\bibitem{miogithub}
Celoria, Daniele.
\newblock PH code
\newblock \url{https://github.com/agnesedaniele/knot-confinement-and-PH}, (2021).

\bibitem{dabrowski2019topoly}
Dabrowski-Tumanski Pawel, Rubach Pawel, Niemyska Wanda, Bartosz~Ambrozy Gren, and Sulkowska Joanna~Ida.
\newblock Topoly: Python package to analyze topology of polymers.
\newblock {\em Briefings in Bioinformatics}, 2019.

\bibitem{de2007coverage}
De Silva Vin and Ghrist Robert
\newblock Coverage in sensor networks via persistent homology
\newblock {\em Algebraic \& Geometric Topology} 7, no. 1 (2007): 339--358.

\bibitem{curvtors}
Diao Yuanan, Claus Ernst, Montemayor Anthony and Ziegler Uta.
\newblock Curvature of random walks and random polygons in confinement
\newblock {\em  Journal of Physics A: Mathematical and Theoretical} 46, no. 28 (2013): 285201.

\bibitem{polyconfinement}
Diao Yuanan, Ernst Claus, Montemayor Anthony and Ziegler Uta.
\newblock Generating equilateral random polygons in confinement
\newblock {\em Journal of Physics A: Mathematical and Theoretical}, 44, no. 40 (2011): 405202.

\bibitem{polyconfinement2}
Diao Yuanan, Claus Ernst, Montemayor Anthony and Ziegler Uta.
\newblock Generating equilateral random polygons in confinement {II}
\newblock {\em Journal of Physics A: Mathematical and Theoretical}, 45, no. 27 (2012): 275203.

\bibitem{polyconfinement3}
Diao Yuanan, Claus Ernst, Montemayor Anthony and Ziegler Uta.
\newblock Generating equilateral random polygons in confinement {III}
\newblock {\em Journal of Physics A: Mathematical and Theoretical}, 45 (2012): 465003.

\bibitem{spectrumconfinement}
Diao Yuanan, Ernst Claus, Montemayor Anthony, Rawdon Eric J. and Ziegler Uta.
\newblock The knot spectrum of confined random equilateral polygons
\newblock {\em Computational and Mathematical Biophysics}, 2(1) (2014).

\bibitem{acnandwrithe}
Diao Yuanan, Ernst Claus, Rawdon Eric J. and Ziegler Uta.
\newblock Average crossing number and writhe of knotted random polygons in confinement
\newblock {\em Reactive and Functional Polymers} 131 (2018): 430-444.

\bibitem{totcurvtors}
Diao Yuanan, Claus Ernst, Rawdon Eric J. and Ziegler Uta.
\newblock Total curvature and total torsion of knotted random polygons in confinement 
\newblock {\em  Journal of Physics A: Mathematical and Theoretical} 51, no. 15 (2018): 154002.

\bibitem{edelsbrunner2010computational}
Edelsbrunner Herbert and Harer John.
\newblock Computational topology: an introduction
\newblock {\em American Mathematical Soc} 4, no. 1 (2010).

\bibitem{stamech}
Flory Paul J. and M. Volkenstein.
\newblock Statistical mechanics of chain molecules
\newblock (1969): 699-700.

\bibitem{ghrist2008barcodes}
Ghrist Robert.
\newblock Barcodes: the persistent topology of data
\newblock {\em Bulletin of the American Mathematical Society} 45, no. 1 (2008): 61--75.

\bibitem{litherland1999thickness}
Litherland Richard~A, Simon Jon, Durumeric Oguz, and Rawdon Eric J.
\newblock Thickness of knots.
\newblock {\em Topology and its Applications}, 91(3):233--244, 1999.

\bibitem{gyrationcompactness}
Lobanov M. Yu, Bogatyreva N. S. and Galzitskaya O. V.
\newblock Radius of gyration as an indicator of protein structure compactness
\newblock {\em Molecular Biology}, 42.4 (2008): 623-628.

\bibitem{biopolymer}
Marenduzzo Davide, Micheletti Cristian and Orlandini Enzo.
\newblock Biopolymer organization upon confinement
\newblock {\em Journal of Physics: Condensed Matter 
}22.28 (2010): 283102.

\bibitem{micheletti}
Micheletti Cristian, Marenduzzo Davide and Orlandini Enzo. 
\newblock Polymers with spatial or topological constraints: Theoretical and computational results
\newblock {\em Physics Reports} 504.1 (2011): 1-73.


\bibitem{ohara}
O'Hara Jun.
\newblock Energy of knots and conformal geometry
\newblock {\em World Scientific, }Vol. 33 (2003).

\bibitem{roadmap}
Otter Nina, Porter Mason A., Tillmann Ulrike, Grindrod Peter and Harrington Heather A.
\newblock A roadmap for the computation of persistent homology
\newblock {\em EPJ Data Science }6 (2017): 1-38.

\bibitem{capsids2}
Petrov Anton S. and Harvey, Stephen C.
\newblock Packaging double-helical DNA into viral capsids: structures, forces, and energetics
\newblock {\em Biophysical Journal} 95, no. 2 (2008): 497-502.

\bibitem{capsids3}
Purohit Prashant K., Kondev, Jan\'e and Phillips Rob.
\newblock Mechanics of DNA packaging in viruses
\newblock {\em Proceedings of the National Academy of Sciences }100, no. 6 (2003): 3173-3178.

\bibitem{rolfsen2003knots}
Rolfsen Dale.
\newblock  Knots and Links
\newblock {\em AMS Chelsea Publishing Series. AMS Chelsea Pub.}, 2003.

\bibitem{knotplot}
Scharein Robert G. and Kellogg S. Booth.
\newblock Interactive knot theory with KnotPlot.
\newblock Multimedia Tools for Communicating Mathematics. Springer, Berlin, Heidelberg, 2002. 277-290.

\bibitem{gyrationselfavoiding}
Shimamura Miyuki K. and Tetsuo Deguchi.
\newblock Finite-size and asymptotic behaviors of the gyration radius of knotted cylindrical self-avoiding polygons
\newblock {\em Physical Review}, E 65.5 (2002): 051802.

\bibitem{stasiak1998ideal}
Stasiak Andrzej and Vsevolod Katritch.
\newblock Ideal knots
\newblock {\em  World Scientific},  Volume~19, 1998.

\bibitem{capsids}
Sun Siyang, Rao Venigalla B., and Rossmann Michael G. 
\newblock Genome packaging in viruses 
\newblock {\em Current Opinion in Structural Biology} 20, no. 1 (2010): 114-120.  

\bibitem{pyknotid}
Taylor Alexander J. and other SPOCK contributors.
\newblock pyknotid knot identification toolkit
\newblock \url{https://github.com/SPOCKnots/pyknotid},  2017.

\bibitem{ripser}
Tralie Christopher, Saul Nathaniel and Bar-On Rann.
\newblock Ripser.py: A lean persistent homology library for python
\newblock {\em Journal of Open Source Software} 3.29 (2018): 925.

\bibitem{gyrationdiffusivity}
Yamamoto Eiji, Takuma Akimoto, Ayori Mitsutake and Ralf Metzler.
\newblock Universal relation between instantaneous diffusivity and radius of gyration of proteins in aqueous solution
\newblock {\em Physical review letters}, 126(12), 128101 (2021).

\bibitem{packing2}
Zinchenko Anatoly, Eisuke Inagaki and Shizuaki Murata.
\newblock Encapsulation of Long Genomic DNA into a Confinement of a Polyelectrolyte Microcapsule: A Single-Molecule Insight
\newblock {\em ACS Omega} 4, no. 1 (2019): 458-464.

\end{thebibliography}
\end{document}